\documentclass[11pt,a4paper]{article}
\usepackage[english]{babel}
\usepackage{comment}
\usepackage[nottoc]{tocbibind}
\usepackage[leqno]{amsmath}
\usepackage{amssymb,mathrsfs,amsthm}
\usepackage{tikz} 
\usepackage{graphicx}
\usepackage{upgreek}
\usepackage{mathtools}
\usepackage{yfonts}
\usepackage{enumitem}
\usepackage{xcolor}
\makeatletter
\newcommand{\leqnomode}{\tagsleft@true\let\veqno\@@leqno}
\newcommand{\reqnomode}{\tagsleft@false\let\veqno\@@eqno}
\newcommand{\bpf}{\begin{proof}}
\newcommand{\epf}{\end{proof}}
\makeatother
\newcommand*{\boldone}{\text{\usefont{U}{bbold}{m}{n}1}}
\newcommand{\w}{\omega}

\newcommand{\eps}{\varepsilon}
\newcommand{\bd}{\mathbf d}
\newcommand{\Gam}{\Gamma}
\newcommand{\vect}{\mathfrak t}
\newcommand{\vecn}{\mathfrak n}
\newcommand{\vecb}{\mathfrak b}
\newcommand{\fancyD}{\mathfrak D}
\newcommand{\fancyE}{\mathfrak E}
\newcommand{\fancyF}{\mathfrak F}
\newcommand{\vbarg}{\overline{v}^g}
\newcommand{\ubarg}{\overline{u}^g}
\newcommand{\chitil}{\widetilde\chi}
\newcommand{\rhotil}{\widetilde\rho}
\newcommand{\etag}{\eta^g}
\newcommand{\wg}{\w^g}
\newcommand{\etace}{\eta^{c}}
\newcommand{\wce}{\w^c}
\newcommand{\etacone}{\eta^{c_1}}
\newcommand{\etactwo}{\eta^{c_2}}
\newcommand{\vg}{v^{g}}
\newcommand{\ug}{u^{g}}
\newcommand{\vcone}{v^{c_1}}
\newcommand{\vctwo}{v^{c_2}}
\newcommand{\ucone}{u^{c_1}}
\newcommand{\uctwo}{u^{c_2}}
\newcommand{\wcone}{\omega^{c_1}}
\newcommand{\wctwo}{\omega^{c_2}}
\newcommand{\nw}{\norm{\w}_{X_T}}
\newcommand{\CC}{\mathcal C}

\newcommand{\EE}{\mathcal E}
\newcommand{\FF}{\mathcal F}
\newcommand{\LTu}{L^2_{uloc}}
\newcommand{\MM}{\mathcal M}
\newcommand{\QQ}{\mathcal Q}

\newcommand{\RR}{\mathbb R}
\newcommand{\NN}{\mathbb N}
\newcommand{\ZZ}{\mathbb Z}
\newcommand{\TT}{\mathbb T}
\newcommand{\RRT}{{\RR^2\times\TT}}
\newcommand{\MMVF}{{\MM^\frac{3}{2}}}

\newcommand{\pa}{\partial}
\newcommand{\na}{\nabla}
\newcommand{\<}{\langle}
\renewcommand{\>}{\rangle}
\DeclareMathOperator{\curl}{curl}
\DeclareMathOperator{\divergence}{div}
\DeclareMathOperator{\Heat}{H}
\DeclareMathOperator{\loc}{loc}
\DeclareMathOperator{\uloc}{uloc}

\DeclareMathOperator*{\pv}{p.v.}

\def\norm#1{\left\vert \left\vert #1  \right\vert \right\vert }

\usepackage{hyperref}
\hypersetup{colorlinks=true, linkcolor = blue,
	anchorcolor = blue,
	citecolor = blue,
	filecolor = blue,
	urlcolor = blue}
 
\newtheorem{lem}{Lemma}[section]
\newtheorem{prop}[lem]{Proposition}
\newtheorem{theo}[lem]{Theorem}
\newtheorem{defi}[lem]{Definition}
\newtheorem{coro}[lem]{Corollary}
\newtheorem{remark}[lem]{Remark}

\usepackage[a4paper]{geometry}
\setlist[enumerate]{label={(\alph*)}}
\numberwithin{equation}{section}

\title{The Cauchy problem for a helical vortex filament in 3D Navier Stokes}
\author{
	Francisco Gancedo
	\\{\footnotesize Departamento de An\'alisis Matemático \& IMUS}
	\\{\footnotesize Universidad de Sevilla}
	\\{\footnotesize Sevilla, Espa\~na}
	\\{\footnotesize email: {\it fgancedo@us.es}}
	\and 
	Antonio Hidalgo-Torné
	\\{\footnotesize Departamento de An\'alisis Matemático \& IMUS}
	\\{\footnotesize Universidad de Sevilla}
	\\{\footnotesize Sevilla, Espa\~na}
	\\{\footnotesize email: {\it ahtorne@us.es}}
}

\begin{document}
\maketitle 

\begin{abstract}
 This paper studies the Cauchy problem for a helical vortex filament evolving by the 3D incompressible Navier-Stokes equations. We prove global-in-time well-posedness and smoothing of solutions with initial vorticity concentrated on a helix. We provide a local-in-time well-posedness result for vortex filaments periodic in one spatial direction, and show that solutions with helical initial data preserve this symmetry. We follow the approach of \cite{bedrossiangermainharropgriffiths23}, where the analogue local-in-time result has been obtained for closed vortex filaments in $\RR^3$. Next, we apply local energy weak solutions theory with a novel estimate for helical functions in non-helical domains to uniquely extend the solutions globally in time. This is the first global-in-time well-posedness result for a vortex filament without size restriction and without vanishing swirl assumptions.
\end{abstract}

{\small
\tableofcontents}

\section{Introduction} 
This paper is devoted to the study of helical vortex filaments in Navier-Stokes equations. 
To this end, we consider the Navier-Stokes equations with one periodic coordinate 
\begin{equation}\label{NSvelocity}
\left\{ 
\begin{aligned}
&\pa_t u(t,x)+u(t,x)\cdot \na u(t,x)-\nu \Delta u(t,x)+\na p(t,x)=0, \quad t\in (0,T), x\in \RRT, \\
&\na \cdot u(t,x)=0, \\
&u(0,x)=u_0(x), 
\end{aligned}
\right.
\end{equation}
where $u_0$ is the prescribed initial data and $\nu> 0$ is the viscosity. We fix $\nu=1$ for the sake of simplicity. Applying the curl operator in the previous system, we obtain the Navier-Stokes equations in vorticity form. Denoting $\w=\na\wedge u$, we get
\begin{equation}\label{NSvorticity}
\left\{ 
\begin{aligned}
&\pa_t \w(t,x)+u(t,x)\cdot \na \w(t,x)-\w(t,x)\cdot\na u(t,x)-\Delta \w(t,x)=0, \quad t\in (0,T), x\in \RRT, \\
&u(t,x)=\na\wedge\Delta^{-1}\w(t,x), \\
&\w(0,x)=\w_0(x)=\na\wedge u_0(x). \\
\end{aligned}
\right.
\end{equation}

In this paper, we consider an initial vorticity $\w_0$ supported on a curve $\Gamma$, which is called a vortex filament. We denote it as 
\begin{equation}\label{vortexfilamentdata}
\w_0=\alpha \delta_{\Gamma},
\end{equation}
where $\alpha\in \RR$ is any constant and the measure $\delta_{\Gamma}$ satisfies 
$$\langle \delta_{\Gamma},\varphi\rangle=\int_\Gamma \varphi\cdot d\vec{s}, \quad \forall \varphi\in C^\infty_c(\RRT).$$

Denoting $\bd(x)$ the distance from $x$ to $\Gamma$, this particular vorticity configuration implies that the velocity field grows as $|\bd(x)|^{-2}$. Therefore, the kinetic energy is infinite near the filament. In the literature, the term vortex filament is often relaxed and refers to vorticities approximately concentrated along a curve.  

The study of vortex filaments in Euler equations, which are \eqref{NSvelocity} with $\nu=0$, dates back to Helmholtz \cite{Helmholtz} and Kelvin \cite{Thomson1868}. In \cite{darios1906} it was shown that an inviscid vortex filament should behave as the binormal flow to leading order. The binormal flow, also known as the self-induction approximation or the vortex-filament equation, has proven to be a very interesting equation by itself. See \cite{jerrardsmets2015,banicavega2020,banicavega2022} and the references therein. In \cite{JerrardSeis2017}, results are given in favor of the vortex filament conjecture. This conjecture states that it is always possible to find solutions of Euler equations whose vorticity remains close to a curve solving the binormal flow. For now, this conjecture is only known to be true in the case of a straight line, a circular line \cite{fraenkel1970}, and a helical line \cite{daviladelpinomussowei2022}.

Regarding the Navier-Stokes equations in $\RR^3$ with initial data including \eqref{vortexfilamentdata}, in \cite{gigamiyakawa1989} it was proved global-in-time existence and uniqueness of smooth solutions for small initial vorticity belonging to the Morrey space $\MMVF$. See Definition \ref{def:morreyR2T}. They consider the vorticity equation in integral form, which can be deduced from \eqref{NSvorticity} using the Duhamel's formula as
\begin{equation}\label{NSintegralvorticity}
\left\{
\begin{aligned}
&\w(t)=e^{t\Delta}\w_0-\int_0^te^{(t-s)\Delta}\divergence (u(s)\otimes \w(s)-\w(s)\otimes u(s))ds, \\
&u=\na\wedge\Delta^{-1}\w, \\
&\w_0(x)=\na\wedge u_0(x).
\end{aligned}
\right.
\end{equation}
Using this formulation, in \cite{bedrossiangermainharropgriffiths23} was proven local-in-time well-posedness of solutions with initial data \eqref{vortexfilamentdata} without restriction on $\alpha$, for any smooth, closed, non-self-intersecting $\Gamma$. The result \cite{kochtataru2001} applies in this setting, but only for small $\alpha$.

Different results that show the existence of solutions with data \eqref{vortexfilamentdata} for arbitrary $\alpha$ need symmetry assumptions that rule out vortex stretching. The case where $\Gamma$ is the vertical line $\{x_1=x_2=0\}$ is known as the Lamb-Oseen vortex, and is essentially a two-dimensional problem. In 2D Navier-Stokes, global-in-time existence with measures as initial vorticity was proven in \cite{cottet1986,giga1988}. Uniqueness for the Lamb-Oseen vortex was obtained in \cite{gallaghergallaylions2005,gallaywayne2005} and later for any measure in \cite{gallaghergallay2005,gallay2012,bedrossianmasmoudi2014}. See also \cite{bedrossiangolding2021} for a 3D approach.
Another case studied is axisymmetric initial data with vanishing swirl. In \cite{fengsverak2015} global-in-time existence and smoothness of solutions was proven, while uniqueness was obtained in \cite{gallaysverak2019}. The short-time dynamics of the unique solution when $\Gamma$ is a circle has been detailed in \cite{gallaysverak2023}, validating the binormal flow approximation. Very recently, the binormal flow approximation has also been validated for any closed non-self-intersecting curve with small $\alpha$ \cite{FontelosVega23}.

An important symmetry conserved by equations \eqref{NSvelocity} is the helical symmetry, which consists of a simultaneous rotation and translation. See Definition \ref{def:helicalsymmetry}. In this case, the natural domain is $\RRT$ or any vertical cylinder centered at the origin. The case $\RRT$ is particularly difficult because in general $u\not\in L^2$ for $\w\in C^\infty_c$, as occurs in 2D.
The cases $u_0\in L^2$ or $H^1$ have been studied in \cite{MahalovTitiLeibovich90}, see also \cite{jiulopesniunussenzveig2018}. They obtain global-in-time well-posedness of weak solutions for $u_0\in L^2$ and strong solutions for $u_0\in H^1$. Using the helicity, one can reduce the system \eqref{NSvelocity} to two dimensions \cite{lopesmazzucatoniunussenzveig2014}, but the resulting system is also complex. In the inviscid case, imposing a vanishing helical swirl \eqref{eq:vanishinghelicalswirl} is a standard assumption in the literature that rules out vortex stretching. In this setting, \cite{dutrifoy1999} showed global-in-time existence of helical solutions for $u_0\in C^s$, for $s>1$. In \cite{ettingertiti2009} global-in-time well-posedness is shown for $\w_0\in L^\infty$. Global-in-time existence was shown in \cite{bronzinussenzveig2015} for vorticity in $L^1\cap L^p$ with $p>4/3$ and in \cite{jiuliniu2017} for $p>1$. Helical solutions solving the vortex filament conjecture for a helical curve were constructed in \cite{daviladelpinomussowei2022}. 
However, the vanishing of the helical swirl is not preserved in the Navier-Stokes equations, which makes it harder to deal with infinite-energy solutions. 

Therefore, in Navier-Stokes with helical symmetry it is natural to consider local energy weak solutions (see Definition \ref{def:localenergyweaksolutions}). Local energy weak solutions were studied in \cite{lemarierieusset2002}, showing global-in-time existence of local energy weak solutions for initial data in $\mathring{E}_2(\RR^3)$ (see Definition \ref{def:uloclpspaces}). This result was extended in \cite{kikuchi2007weak}, adding a force term. The global-in-time existence of local energy weak solutions have been obtained in the half space in \cite{maekawamiuraprange2019}. In \cite{basson2006}, global-in-time well-posedness and uniqueness were proven for data in $\mathring{E}_2(\RR^2)$. The study of the pressure requires special attention in local energy weak solutions, see \cite{fernandezdalgolemarierieusset2021,bradshawtsai2022}. For other considerations and more results of existence for data in $\LTu(\RR^3)$, see \cite{bradshawtsai2020} and the references therein.

In this work, we prove the first global-in-time well-posedness result for a vortex filament without vanishing swirl assumptions. We achieve this for $\w_0$ given by $\eqref{vortexfilamentdata}$ with helical $\Gamma$.
First, we prove local-in-time well-posedness for any $\Gamma\subset\RRT$ smooth, $x_3$-periodic, non-self-intersecting. The Biot-Savart and other kernels in $\RRT$ exhibit very different behavior at infinity than in $\RR^3$ \cite{bronzinussenzveig2015,jiuliniu2017}. We overcome the difficulties and adapt the proof of \cite{bedrossiangermainharropgriffiths23}, where a local-in-time well-posedness result is proven for $\Gamma\subset\RR^3$. Another nontrivial question is the preservation of helical symmetry. We show this property with a careful study of the symmetry of each component in which the vorticity is decomposed.

The velocity field obtained for short times belongs to $\mathring{E}_2(\RRT)$ for $t>0$. We show that helical divergence-free $L^2$ vector fields are dense in helical $\mathring{E}_2(\RRT)$ fields. We also prove a new helical estimate in non-helical domains. With these, we adapt the construction of \cite{kikuchi2007weak} to obtain a global-in-time local energy weak solution, and then perform energy estimates using helical estimates to obtain uniqueness. Finally, we prove that the local-in-time solution obtained for a helical vortex filament can be uniquely continued globally in time, remaining smooth. In subsection \ref{subsection:statementandoutline} we sketch the proof and its difficulties.

\subsection{General notation}\label{subsection:generalnotation}
In this subsection we clarify some notation. Given two comparable quantities $a,b$, we denote $a\lesssim b$ if there exists a constant $C>0$ such that $a\leq Cb$. 

We will use $x=(x_1,x_2,x_3)$ or $y=(y_1,y_2,y_3)$ as cartesian coordinates, and $(\rho,\theta,z)$ as cylindrical coordinates. We denote by $\widetilde{x}$ the first two components of the vector $x$. 

We define the bump function $\chi_l=\chi_l(|\widetilde{x}|)$ as a smooth, non-negative bump function identically equal to $1$ on the set $\{|\widetilde{x}|\leq l\}$ and supported on the set $\{|\widetilde{x}|\leq 2l\}$. 

We also use the Japanese bracket $\langle x\rangle=(1+|x|^2)^{1/2}.$

In this paper we will only use the Fourier transform in the $x_3$ variable. We will denote $\hat{u}(x_1,x_2,\xi)$ the Fourier transform of $u(x_1,x_2,x_3)$ in its third variable.

We denote by $B(x,l)$ the ball with center $x$ and radius $l$. This will be used for balls in $\RR^3$ and $\RRT$. In $\RRT$ we also work with cylinders. If $x=(x_1,x_2,x_3)\in \RRT$, we denote by $C(x,l)$ the cylinder $C(x,l)\coloneqq \{(y_1,y_2,y_3)\in \RRT : |(x_1,x_2)-(y_1,y_2)|<l\}.$ 

We define the bilinear operator $B[f,g]=\divergence(f\otimes g -g\otimes f)$.

\subsection{Statement of results and outline of the proof}\label{subsection:statementandoutline}
In this subsection, we will give some definitions needed to state the main results and explain the ideas and difficulties of the proofs. Other less general definitions will be given when they are needed.

We adapt the usual definition of Morrey spaces to the domain $\RRT$ using \cite{geisler88}. 
\begin{defi}\label{def:morreyR2T}
Let $1\leq q\leq p$. $M^p_q$ is the space of functions $f\in L^p_{\loc}$ such that 
$$\norm{f}_{M^p_q}\coloneqq \sup_{y\in \RRT, l>0}|B(y,l)\cap (\RRT)|^{\frac{1}{p}-\frac{1}{q}}\left(\int_{|B(y,l)\cap (\RRT)|}|f|^qdz\right)^\frac{1}{q}<\infty,$$
and $\MM^p$ the space of Borel measures $\mu$ such that 

\begin{equation}\label{eq:Morreyspacedef}
\norm{\mu}_{\MM^p}\coloneqq\sup_{y\in \RRT, l>0}|B(y,l)\cap (\RRT)|^{\frac{1}{p}-1}\left|\mu(B(y,l)\cap (\RRT)) \right|<\infty.
\end{equation}
\end{defi}

Note that a vortex filament \ref{vortexfilamentdata} is embedded in $\w\in\MM^\frac{3}{2}$. We define a solution to a vortex filament initial data as in \cite{bedrossiangermainharropgriffiths23}.
\begin{defi}\label{def:mildsolution}
Let $T>0$. We say that $\w$ is a mild solution to the Navier-Stokes equations with initial data $\w_0\in\MMVF$ if
\begin{enumerate}
    \item $\w\in \CC_w([0,T],\MMVF)$, with $\w(0)=\w_0$.
    \item $\w$ satisfies equation \eqref{NSintegralvorticity}
    \item $\w$ is divergence free in the sense of distribution for all $t\in [0,T]$. 
\end{enumerate}
\end{defi}
Our first result is the analogue in $\RRT$ of the main result in \cite{bedrossiangermainharropgriffiths23}, which is done in $\RR^3$.
\begin{theo}[Local-in-time existence $x_3$-periodic vortex filaments]\label{theo:localexistence}
For any $\alpha\in \RR$, and any $x_3$-periodic, smooth, and non-self-intersecting curve $\Gam$ there exist $T>0$ and a mild solution to the Navier-Stokes equations in the sense of definition \ref{def:mildsolution} with initial data \eqref{vortexfilamentdata}. The solution $\w\in C^0((0,T],L^1\cap L^2)$, with $u \in C^0((0,T],L^p)$, for $2<p<6$.
\end{theo}
\begin{remark}
The initial data can be taken more general. As in \cite{bedrossiangermainharropgriffiths23}, a noise can be added to the initial data and the result still holds. However, as noted in \cite{bedrossiangermainharropgriffiths23}, the addition of noise is not the most challenging part. The idea of the proof is to construct a contractive fixed point scheme, so we also have a uniqueness result as a corollary. See corollary \ref{coro:uniquenesssmalltime}.
\end{remark}
We adapt the strategy in \cite{bedrossiangermainharropgriffiths23}, which is to use a fixed-point scheme. Since the initial data is not small, it is needed a good approximate solution $\w^g$ so that the remainder $\w-\w^g$ is small. Near the curve, the solution should be similar to the solution of the straight filament $\eta^g$, known as the Lamb-Oseen vortex, because the curvature is expected to be subcritical. The first step is to construct the approximate solution $\w^g$ using the Lamb-Oseen vortex and a good mapping between the straight filament and the curve $\Gamma.$ We also need control of the linearized equation around the Lamb-Oseen vortex, which controls the evolution for short times near the filament. The distinction of near and far from the filament is made rigorous by splitting the unknown vorticity in two, and treating each part appropriately. 

The most difficult part in the adaptation of the proof is that the kernels of the heat and Biot-Savart operator are no longer explicit in $\RRT$ and also have different decay at infinity, being similar to the decay in two dimensions. This changes the integrability of the velocity while keeping the same integrability in the vorticity. For example, in \cite{bedrossiangermainharropgriffiths23} the velocity field belongs to $L^2$ for any positive time, but in our case it does not for any time. Note that in $\RRT$ the Biot-Savart law decays as $|\widetilde{x}|^{-1}$, which is not integrable at infinity in $\RRT$, as occurs in 2D. 

In Section \ref{Section:constructionfixedpoint} we construct the mapping that transforms vorticities and velocities and define the fixed point scheme. In Section \ref{section:estimatesandlocalexistence} we perform estimates to show that the fixed point is contractive, and then we prove Theorem \ref{theo:localexistence}. 

In the case where $\Gamma$ has helical symmetry, we get additional results. We define the helical symmetry following \cite{bronzinussenzveig2015}. Helical symmetry implies periodicity in $x_3$ with period $\kappa$. We fix $\kappa=2\pi$ for simplicity. 
\begin{defi}\label{def:helicalsymmetry}
Denote as $R_{\theta}$ the rotation matrix around the $x_3$ axis
with angle $\theta$
\begin{gather*}
 R_\theta=
 \begin{pmatrix} \cos\theta & \sin\theta & 0 \\ -\sin\theta &\cos\theta &0 \\
 0 & 0 & 1
 \end{pmatrix},
\end{gather*}
and $S_\theta$ the superposition of $R_\theta$ and a vertical translation of size $\theta$
\begin{gather}\label{eq:helicalStheta}
 S_\theta(x)=
 R_\theta(x)
 +\begin{pmatrix} 0 \\ 0 \\
 \theta
 \end{pmatrix}.
\end{gather}
We say that a scalar function $f(x):\RR^3\mapsto \RR$ is helical if 
$$f(x)=f(S_\theta x)\quad \forall x\in \RR^3, \theta\in \RR.$$
We say that a vector function $v(x):\RR^3\mapsto\RR^3$ is helical if 
$$R_\theta(v(x))=v(S_\theta x) \quad \forall x\in \RR^3, \theta\in \RR.$$
We say that a non-smooth function is helical if it is limit of helical smooth functions in the sense of distributions. The same applies to vector functions. 
We say that a helical vector field $v(x)$ has vanishing helical swirl if 
\begin{equation}\label{eq:vanishinghelicalswirl}
v(x)\cdot (x_2,-x_1,1) =0.
\end{equation}
\end{defi}

\begin{remark}
Using cylindrical coordinates $(\rho, \theta, z)$ and making the change $\theta\to \theta+z=\xi$, the previous definition is equivalent to saying that a scalar function $f$ is helical when $f(\rho,\xi,z)$ does not depend on the value of $z$. 
\end{remark}

Although Navier-Stokes equations preserve helical symmetry, if $\w_0$ in \eqref{vortexfilamentdata} is helical, we cannot directly deduce that the solution obtained in Theorem \ref{theo:localexistence} is helical because the initial data is rough. We can prove this in a proposition.

\begin{prop}\label{prop:helicalsymmetrypreserved}
Let $\alpha\in \RR$, $\Gamma$ be a helical curve, and $\w_0$ given by \eqref{vortexfilamentdata}. Then, the solution obtained in Theorem \ref{theo:localexistence} is helical during its time of existence.
\end{prop}
One possible strategy could be to mollify the initial data to approach the solution. It is clear that this strategy does not work because the norms used in the fixed-point theorem do not behave well with mollification because of the irregularity of the solution at time $0$. The challenge in this proposition is that we have split the vorticity into three parts, and that we work in two different frames simultaneously. To show the preservation of symmetry, we prove that being helical in the physical frame is equivalent to being independent of $x_3$ in the straightened frame. The fixed-point scheme itself can be seen as an iterative process, and we can show that each iteration preserves the symmetry. Starting the process with $\w^g$, we conclude that the solution is also symmetric. This study is carried out in Section \ref{Section:preservationsymmetry}. The remaining theorems are proven in Section \ref{section:globaltimeexistence}.

We can uniquely extend the solution obtained when $\Gamma$ is helical using local energy solutions. We define the loc-uniform Lebesgue spaces as in \cite{lemarierieusset2002},

\begin{defi}\label{def:uloclpspaces}
For any $1\leq p \leq \infty$,
$$L^p_{\uloc}(\RRT)=\{u\in L^2_{\loc} : \norm{u}_{L^p_{\uloc}}(\RRT)=\sup_{x_0\in\RRT}\norm{u}_{L^p(B(x_0,1))}<\infty \}.$$
We define also $\mathring{E}_p$ as the closure of divergence free compactly supported smooth functions with respect to the norm $L^p_{\uloc}$.     
\end{defi}
In \cite{lemarierieusset2002} it is shown that $\mathring{E}_p$ coincides with the functions in $L^p_{\uloc}$ that tend to $0$ at infinity. Since the initial data does not have finite energy, Leray-Hopf solutions does not make sense and we need another definition. We follow \cite{kikuchi2007weak} and \cite{bradshawtsai2022}.

\begin{defi}[Local energy weak solutions]\label{def:localenergyweaksolutions}
Fix $0<T<\infty$. A vector field $u\in L^2_{\loc}([0,T]\times \RRT )$ is a local energy weak solution to \eqref{NSvelocity} with divergence free initial data $u_0\in \LTu(\RRT)$ if:
\begin{enumerate}
    \item $u\in L^\infty([0,T]; \LTu)$, $\na u\in \LTu([0,T]\times \RRT )$.
    \item For some $p\in L^{3/2}_{\loc}([0,T]\times \RRT)$, the pair $(u,p)$ is a distributional solution to \eqref{NSvelocity}.
    \item For any compact $K\subset \RRT$, $u(t)\to u_0$ in $L^2(K)$ as $t\to 0^+$. 
    \item The function $t\mapsto \int_\RRT u(x,t)\cdot h(x)dx$ is continuous in $t\in [0,T)$ for any compactly supported $h\in L^2(\RRT)$.
     \item For almost any $0<t<T$ and all nonnegative $\phi\in C^\infty_c((0,T)\times \RRT )$,
     $$\int_\RRT |u(t)|^2\phi(t) dx+2\int_0^t\int_\RRT|\na u|^2\phi dx ds
     \leq \int_0^t\int_\RRT|u|^2(\pa_s\phi+\Delta \phi)+(|u|^2+2p)(u\cdot\na)\phi dx ds.$$
     \item For any $x_0\in \RRT$ there exist a function $c_{x_0}\in L^{3/2}(0,T)$ such that for any $x\in C(x_0,3/2)$ and $0<t<T$,
     \begin{equation}\label{eq:pressuredecomposition}
     \begin{aligned}
    p(t,x)-c_{x_0}=&-\frac{1}{3}|u(t,x)|^2+\pv \int_{C(x_0,2)}N_{i,j}(x-y)(u_i(s,y)u_j(s,y))dy \\
    &+\int_{(\RRT)\setminus C(x_0,2)}(N_{i,j}(x-y)-N_{i,j}(x_0-y))(u_i(s,y)u_j(s,y))dy,
    \end{aligned}
     \end{equation}
      where $N$ is the Newtonian potential and $N_{i,j}=\pa_i\pa_jN$. 
\end{enumerate}
\end{defi}
\begin{remark}
In the last property of the previous definition we use cylinders instead of balls because we work in $\RRT$, but this change is not relevant. The radius of the cylinder in the pressure decomposition is not important and can be taken arbitrarily big, we fix it for ease of reading. The pressure decomposition formula can be deduced from the other conditions; see \cite{kikuchi2007weak} or \cite{bradshawtsai2022} for a more general discussion.
\end{remark}

\begin{theo}\label{theo:globalexistenceintroduction}
Let $T>0$, $u_0\in \mathring{E}_2$ helical. Then, there exists a unique helical local energy weak solution to \eqref{NSvelocity}. 
\end{theo}
In \cite{basson2006} it is used that there is no vortex stretching in two dimensions, so we need another approach. In \cite{kikuchi2007weak} local-in-time existence of local energy weak solutions is proven by regularizing the initial data and taking limits. Then, it is shown that the solution can be continued for all times. To use this result, we need to regularize the initial data while preserving the symmetry. This is done in Lemma \ref{lem:approxlocalenergy}. Next, we perform local energy estimates to show uniqueness. For that, we need to use a novel helical interpolation \eqref{eq:helicalinterpolation} in non-helical domains. 

To conclude the paper, we apply energy estimates to the local-in-time solution obtained for a helical vortex filament to prove the following theorem.

\begin{theo}[Global-in-time well-posedness for a helical vortex filament]\label{theo:globalwell-posednessintroduction}
For any $\alpha\in\RR$ and any $\Gamma$ helical curve, there exists a unique global-in-time smooth helical solution to the Navier-Stokes equations with initial vorticity given by \eqref{vortexfilamentdata}.
\end{theo}

\section{Construction of fixed point scheme}\label{Section:constructionfixedpoint}
The purpose of this section is to construct an approximate solution to the problem and then to set the fixed-point scheme. The content is analogue to \cite[Section 6]{bedrossiangermainharropgriffiths23}
with the difference that we work on $\RRT$ instead of $\RR^3$. Our curve is z-periodic instead of closed with length $2\pi$, but the analysis needed is local and can be applied in our case. We will follow a similar development and notation and skip the proofs, which can be found in \cite{bedrossiangermainharropgriffiths23}. 

\subsection{Geometrical setting}\label{subsection:geometricalsetting}
The purpose of this subsection is to prepare the ingredients needed to construct the approximate solution. 
 Given a smooth z-periodic curve $\Gam$, we can find a constant speed parameterization $\gamma:\TT\to \RR^2\times\TT$. 
The constant speed $c_s\coloneqq\left|\pa_s \gamma(s)\right|$ will be proportional to the length of the curve over one period. 
Consider the canonical orthonormal frame along $\Gam$, that is, $\vect, \vecn, \vecb :\TT\to \RR^3$ such that $c_s\vect=\gamma', \vecn\propto \vect'$ and $\vecb=\vect\times\vecn.$  

For $0<r\ll 1$ to be fixed later, we define a neighborhood $\Gam_r$ around the curve and a straightened frame $\Sigma_r$ as $$\Gam_r\coloneqq\{y\in \RR^2\times\TT:dist(y,\Gam)<32r\},$$
$$\Sigma_r\coloneqq\{x\in \RR^2\times\TT:|\widetilde{x}|<32r\}.$$ 
Then define $\Phi:\Sigma_r\to\Gam_r$ as
 $$\Phi(x)\coloneqq\gamma(x_3)+x_1\vecn(x_3)+x_2\vecb(x_3) ,
 $$
which will be bijective if $r$ is small enough. From now on, we will assume without saying that $r$ is small enough so that $\Phi$ is bijective. Taking derivative respect to $x_3$ in the previous equation, we obtain

$$\pa_z\Phi(x,z)=c_s\vect+x_1(\vecn'\cdot \vect)\vect+x_1(\vecn'\cdot \vecb)\vecb+x_2(\vecb'\cdot \vect)\vect+x_2(\vecb'\cdot \vecn)\vecn\eqqcolon\fancyD \vect +\fancyE\vecn+\fancyF\vecb,$$
so the jacobian of $\Phi$ is 

$$J=\nabla\Phi=\left[ \vecn \quad \vecb \quad \fancyD \vect +\fancyE\vecn+\fancyF\vecb\right].$$
The determinant of $J$, due to orthogonality, is $\fancyD\approx c_s$ in the domain.
Now, we define the maps that will be used to transform velocity fields ($P_\Phi$) and vorticity fields ($Q_\Phi$) from the straightened frame to the curved frame:
\begin{equation}\label{eq:defPQ}
	\begin{aligned}
		P_\Phi v(y) &= Jv(\Phi^{-1}(y)), \\
		Q_\Phi \eta(y) &= \fancyD^{-1}J\eta(\Phi^{-1}(y)). \\
	\end{aligned}
\end{equation}
The reason for using different mappings is that they have good properties together, as we can see in the following lemma.
\begin{lem}\label{lemma:operatorstransformations}
	Let $v,\eta:\Sigma_{r}\to \RR^3$ be smooth vector fields defined in the straightened frame. We have the identities:
\begin{enumerate}
\item \underline{Divergence operator.} 
	\begin{equation*}
		\nabla \cdot (Q_\Phi\eta)=(\fancyD^{-1}\nabla\cdot\eta)\circ\Phi^{-1}.
	\end{equation*}
\item \underline{Curl operator.}
	\begin{equation*}
	\nabla \times (Q_\Phi\eta)=Q_\Phi \curl_\Phi\eta,
	\end{equation*}
where the modified curl operator is defined as 
	\begin{equation*}
	\curl_\Phi\eta\coloneqq\na\times\eta+E^j\pa_j\eta+F\eta,
	\end{equation*}
where the matrix $F$ is smooth and bounded, and the matrices $E^j$ satisfies bounds that we will state in the next lemma.
\item \underline{Bilinear operator.}
	\begin{equation*}
	B[P_\Phi v,Q_\Phi\eta]=Q_\Phi B[v,\eta].
	\end{equation*}
\item \underline{Laplacian.}
	\begin{equation*}
	\Delta Q_\Phi\eta=Q_\Phi\Delta_\Phi\eta,
	\end{equation*}
where the modified Laplacian is defined as
	\begin{equation*}
	\Delta_\Phi\coloneqq\Delta+A^{i,j}\pa_i\pa_j+B^j\pa_j+C,
	\end{equation*}
where the matrices $B^j, C$ are smooth and bounded and the matrix $A$ satisfies bounds that we will state in the next lemma.
\end{enumerate}
\end{lem}
We detail now 
the estimates of different elements from this subsection that will be used later. 

\begin{lem}\label{changevariablebounds}
	Provided $r$ is small enough and $x\in \Gam_{r}$, for any $j\geq0$, $i\in \{1,2,3\}$ and multi-index $\alpha\in \NN^3$,
	\begin{equation*}
		\begin{aligned}
		&	|\fancyD|\gtrsim 1, \\
		&	|\pa^j_{x_3}(\fancyD-c_s)|+|\pa_{x_3}^jE^i|+|\pa_{x_3}^jA|\lesssim|\widetilde{x}|,\\
		& |\na^\alpha_{x}\Phi|+|\na^\alpha_{x}X|+
		   |\na^\alpha_{x}\fancyD|\lesssim|\widetilde{x}|,
		\end{aligned}
	\end{equation*}
where $X$ is any of $J, A, B^i, C, E^i, F$.
\end{lem}
Notice that due to the previous bounds the operators $\curl$ and $\curl_\Phi$ are equal at higher order, except for terms that are small on $\Gam_{r}$. The same occurs for $\Delta$ and $\Delta_\Phi$. The last result of this subsection is consequence of the previous lemma.

\begin{lem}
 Provided $r$ is small enough, for any $x,x'\in \Sigma_{r}$ we have 
 $$|\Phi(x)-\Phi(x')|\approx |\widetilde{x}-\widetilde{x}'|+|x_3-x_3'|.$$
\end{lem}
We will simultaneously work on the straightened frame and in the curved frame. The previous lemma together with lemma \ref{changevariablebounds} will help us to control the change of variables.

Since the change of variables is only available close to the curve $\Gam$, we need a bump function to control it rigorously. In the straightened frame we use the bump function $\chi_r$ defined in subsection \ref{subsection:generalnotation}. The corresponding bump functions on the curved frame $\widetilde\chi_r$ are defined as $\widetilde\chi_r\circ\Phi \coloneqq \chi_r$. Observe that with this definition, we have 
$$\widetilde\chi_rP_\Phi( \cdot ) = P_\Phi(\chi_r \cdot),\quad \widetilde\chi_rQ_\Phi( \cdot ) = Q_\Phi(\chi_r \cdot).$$

\subsection{Approximate solution}
First of all, since we are going to work with vorticities and velocities defined on the straightened and curved frame, we fix more notation. We will denote with $\eta$ vorticities defined on the straightened frame, and $\w$ vorticities on the curved frame. For velocities, we will use $v$ for the ones defined on the straightened frame, and $u$ for the ones on the curved frame. 

To make the fixed point argument work, we need to construct a good initial approximation so the remaining terms are small. As we remarked in subsection \ref{subsection:statementandoutline}, a suitable mapping of the Lamb-Oseen vortex might be a good candidate. We denote the Lamb-Oseen vortex and its correspondent velocity as
\begin{equation}\label{defetagvg}
	\begin{aligned}
	\eta^g(t,x)&=\left[0 \enspace 0 \enspace \frac{\alpha}{4\pi t}e^{-\frac{|\widetilde{x}|^2}{4t}}\right]^t, \\
	\vbarg(t,x)&=\left[\frac{\alpha \widetilde{x}^\bot}{2\pi|\widetilde{x}|^2}\! \left(1-e^{-\frac{|\widetilde{x}|^2}{4t}}\right) \enspace 0 \right]^t,
	\end{aligned}
\end{equation}
where $(x_1,x_2)^\bot\coloneqq (-x_2,x_1)$. The corresponding magnitudes on the curved frame are denoted as
\begin{equation}
\label{def:wgubarg}
\wg=Q_\Phi(\chi_{2r}\etag), \quad \ubarg=P_\Phi(\chi_{4r}\vbarg).
\end{equation}
We must check which equation satisfies $\wg$ to know which equation must satisfy the remaining terms. Note that $\etag$ solves the heat equation since the bilinear term vanishes due to symmetry. Mollification commutes with derivatives, and again by symmetry the bilinear term vanishes for $\etag$, so we can write
$$\pa_t\etag+B[\chi_r\vbarg,\etag]=\Delta\etag$$
and now apply $Q_\Phi(\chi_{2r}\cdot)$ and Lemma \ref{lemma:operatorstransformations} to obtain

\begin{equation}\label{eq:evwg}
	\pa_t\wg+B[\chitil_r\ubarg,\wg]=\Delta\wg+\underbrace{Q_\Phi(\chi_{2r}\Delta \etag)-\Delta\wg}_{\EE^g}.
\end{equation}

\subsection{Unknowns and fixed point scheme}

By construction, the vorticity $\wg$ attains as initial data $\delta_\Gam$, so we can expect the remaining vorticity $\wce$ to be small for short times provided that equation \eqref{eq:evwg} is a good approximation to NS. It is convenient to split the remaining vorticity into two terms: $\wcone$ to track the evolution close to the curve, and $\wctwo$ to do the rest. We will keep the same superindexes for the associated vorticity in the straightened frame and velocities. We leave without superindex the total vorticity $\omega=\wg+\wce=\wg+\wcone+\wctwo$. Regarding velocities, $u$ will always be related to the corresponding $\w$ by the Biot-Savart law $$u^*=\nabla\times(-\Delta)^{-1}\w^*,$$ and except for $\vbarg$ and $\ubarg$, $v$ and $u$ will always be related by 
\begin{equation}\label{eq:relationuv}
 v^*\coloneqq P_\Phi^{-1}(\chitil_ru^*).   
\end{equation}
With respect to vorticities, we choose the following relations
\begin{equation}\label{eq:wetarelations}
\wg\coloneqq Q_\Phi(\chi_{2r}\etag), \quad 
\wcone \coloneqq Q_\Phi(\chi_{2r}\etacone), \quad 
\etactwo \coloneqq Q_\Phi^{-1}(\chitil_{2r}\wctwo).
\end{equation}

We define $\etacone$ through the equation

\begin{equation}\label{eqevetarone}
\pa_t\etacone+B[v,\etace]+B[v-\chi_r\vbarg,\etag]=\Delta\etacone,
\end{equation}
with null initial data. This implies that (applying $Q_\Phi(\chi_{2r}\cdot)$)

\begin{equation}\label{eqevwrone}
\pa_t\wcone+B[\chitil_ru,\wce]+B[\chitil_r(u-\ubarg),\wg]=\Delta\wcone+\underbrace{Q_\Phi(\chi_{2r}\Delta \etacone)-\Delta\wcone}_{\EE^{c}}.
\end{equation}
The remaining $\wctwo$ must close Navier-Stokes on the curved frame, which is the physical one, so
\begin{equation}\label{eqevwrtwo}
	\pa_t\wctwo+B[(1-\chitil_r)u,\w]=\Delta\wctwo-\EE^{g}-\EE^{c}.
\end{equation}
To set the fixed point scheme, the next step is to define the norm for each variable. 

\subsubsection{Norms}
We will use the same norms as \cite{bedrossiangermainharropgriffiths23}, with the difference that for us $y\in \RR^2\times \TT$ instead of $\RR^3$. We will have one norm $	\norm{\cdot}_{N^\beta}$ defined on the straightened space, to be used with terms close to $\Gam$, and another norm $\norm{\cdot}_F$ for the physical space. Their definitions are 

\begin{equation}\label{eq:normsdefinition}
	\begin{aligned}
	\norm{\eta^*}_{N^\beta} &\coloneqq\sup_{0<t\leq T} \sqrt{t}\norm{\left<\frac{\widetilde{x}}{\sqrt{t}}\right>^m\<\sqrt{t}\nabla\>^\beta\eta^*}_{B_{x_3}L^2_x}, \\
	\norm{\omega^*}_F &\coloneqq\sup_{0<t\leq T} \sqrt{t}\norm{(1-\chitil_{6r})\left<\frac{\bd}{\sqrt{t}}\right>^m\omega^*}_{L^3_y}, \\
		\end{aligned}
\end{equation}
where $\beta\in (0,\frac{1}{4})$ is a real number, $m\geq 2$ is a fixed integer, the function $\bd(y)\coloneqq \text{dist}(y,\Gam),$ and the norm 
$$\norm{f}_{B_{x_3}L^p_x}\coloneqq \sum_{\zeta\in\ZZ}\norm{\hat{f}(x_1,x_2,\zeta)}_{L^p_{x_1,x_2}}.$$
We are making a slight abuse of notation, the above norm should be written as $B_{x_3}L^p_{\widetilde{x}}$ or $B_{x_3}L^p_{x_1,x_2}$, but there is no confusion since we are also applying a norm in $x_3$.
For fractional $\beta$, we understand $\na^\beta$ or $\langle \sqrt{t}\na\rangle^\beta$ as a Fourier multiplier.

The norm that will be used at the fixed point argument is

\begin{equation} \label{normdefinition}  
\norm{\w}_{X_T}\equiv\norm{(\etacone,\wctwo)}_{X_T}\coloneqq\norm{\etacone}_{N^\beta}+M\norm{Q_\Phi^{-1}[\chitil_{8r}\wctwo]}_{N^0}+M\norm{\wctwo}_F,
\end{equation}
where $M$ is a constant that we will choose later. We will prove that Banach's fixed point theorem works in the ball 

$$B_{\eps,T,r,M}\coloneqq\left\{(\etacone,\wctwo):\norm{(\etacone,\wctwo)}_{X_T}\leq \eps \right\}$$
if constants are chosen appropiately.

\subsubsection{Mapping from Fixed Point Theorem}
Now, we derive the fixed-point mapping from equations \eqref{eqevetarone}, \eqref{eqevwrtwo}. We will use Duhamel's formula. For that, we choose an evolution operator for each variable, then rewrite equations \eqref{eqevetarone}, \eqref{eqevwrtwo} coherently. 

Let $S(t,s)$ be the solution operator associated to the Navier-Stokes linearized equation around the Lamb-Oseen vortex
\begin{equation}
\eta(t)=S(t,s)\eta(s), \quad \left\{
		\begin{aligned}
		&\pa_t\eta+B[\vbarg,\eta]+B[\textbf{v},\eta^g]=\Delta \eta \\
	&	\textbf{v}=-\nabla\times(\Delta)^{-1}\eta.
		\end{aligned}
	\right.
\end{equation}
We have used $\textbf{v}$ instead of $v$ in order to avoid confusion, recall that $v^*\neq-\nabla\times(\Delta)^{-1}\eta^*$, Biot-Savart law only holds for $u^*$ and $\w^*$. We can rewrite \eqref{eqevetarone} as 

\begin{equation*}
	\begin{aligned}
\pa_t\etacone &+B[\vbarg,\etacone]+B[(-\Delta)^{-1}\nabla\times\etacone,\eta^g]=\Delta\etacone\\
&-B[v-\vbarg,\etacone]-B[v,\etactwo]
-B[\vg-\chi_R\vbarg,\etag]-B[v-\vg-(-\Delta)^{-1}\nabla\times\etacone,\etag] 
	\end{aligned}
\end{equation*}

The first part of our mapping will be 
\begin{equation}\label{eq:a1}
a^{1}=-\int_0^tS(t,s)f^{c_1}(s)ds,
\end{equation}
where
\begin{equation*}
	\begin{aligned}
		f^{c_1}=B[v-\vbarg,\etacone]+B[v,\etactwo]
		+B[\vg-\chi_r\vbarg,\etag]+B[v-\vg-(-\Delta)^{-1}\nabla\times\etacone,\etag] 
	\end{aligned}
\end{equation*}

Denoting by $e^{t \Delta}$ the heat operator, we set from \eqref{eqevwrtwo} the second part of our mapping as

\begin{equation}\label{eq:a2}
a^{2}=-\int_0^te^{(t-s)\Delta}f^{c_2}(s)ds,
\end{equation}
where
$$
f^{c_2}=B[(1-\chitil_r)u,\w]+\EE^{g}+\EE^{c}.
$$

By construction, a fixed point of the mapping $\QQ(\etacone,\wctwo)\coloneqq(a^{1},a^{2})$ will solve \eqref{eqevetarone} and \eqref{eqevwrtwo}.

\section{Fixed Point Estimates, Short Time Existence}\label{section:estimatesandlocalexistence}
In this section we first prove the estimates needed to make the fixed-point argument work. We will control several quantities with the norm $\norm{\cdot}_{X_T}$ and constants $r, T, M$. We assume that $r<1, M>1$. Any other constant will be implicitly included in the notation $\lesssim$. Then, we will conclude the local-in-time existence of solutions. We adapt \cite[Section 7]{bedrossiangermainharropgriffiths23}. Since we work on a different domain, the kernels behave differently and some estimates change. For ease of reading, we try to include at least a short proof of each estimate.
We denote the cartesian coordinates $x$ in the straightened frame and $y$ in the physical frame.

\subsection{Preliminary estimates} 
In this subsection we show estimates that are derived directly from definitions of norms and spatial estimates, i.e., we do not use the evolution equations. We do not intend to make them sharp, we just present what we will use later.

\begin{lem}[Weight bounds]\label{lem:weightbounds}\leavevmode 
Let $0<c$, $1< p\leq \infty$.
\begin{align}
\norm{(1-\chi_{cr})\<t^{-\frac{1}{2}}\widetilde{x}\>^{-m}}_{B_{x_3}L^p_x}+\norm{(1-\chitil_{cr})\<t^{-\frac{1}{2}}\bd\>^{-m}}_{L^p_y}&\lesssim t^{m/2}r^{2/p-m}\label{estweightfar} \\
\norm{\<t^{-\frac{1}{2}}\widetilde{x}\>^{-m}}_{B_{x_3}L^p_x}&\lesssim t^{1/p} \label{estweightfull}
\end{align}
\end{lem}
\bpf
For \eqref{estweightfar} use that far from $0$, $\< \cdot \>$ behaves like $|\cdot|$. For \eqref{estweightfull} use the change of variables $x'=xt^{-1/2}$. 
\epf 

\begin{lem}[Vorticity bounds in straightened frame]\label{lem:boundsstraightenedframe}\leavevmode 
	
\begin{enumerate}
\item \underline{Full space bounds}	
\begin{align}
\norm{\etag}_{N^\beta}
+t^\frac{1}{4}\norm{\etag}_{B_{x_3}L_x^{4/3}}
+\norm{\< \sqrt{t}\nabla \> ^\beta \etag}_{B_{x_3}L^1_x}
&\lesssim 1 \label{estetagBz} \\
\norm{\< \sqrt{t}\nabla \> ^\beta\etacone}_{B_{x_3}L^1_x}+
t^\frac{1}{4}\norm{\etacone}_{B_{x_3}L_x^{4/3}}
 &\lesssim \nw \label{estetar1Bz}\\
\norm{Q_\Phi^{-1}(\chitil_{8r}\wctwo)}_{B_{x_3}L^1_x}+
t^\frac{1}{4}\norm{\etactwo}_{B_{x_3}L_x^{4/3}}
+t^\frac{1}{4}\norm{Q_\Phi^{-1}(\chitil_{8r}\wctwo)}_{B_{x_3}L_x^{4/3}}&\lesssim M^{-1}\nw \label{estetar2Bz}\\
t^\frac{1}{2}\norm{\<t^{-\frac{1}{2}}\widetilde{x}\>^m \etacone}_{L^2_{x}}
+M^{-1}t^\frac{1}{2}\norm{\<t^{-\frac{1}{2}}\widetilde{x}\>^m \etactwo}_{L^2_{x}}&\lesssim \nw \label{estetar1r2L2w}
\end{align}
 
\item \label{estetaBzL43}\underline{Bounds far from filament}
For $c>0$ and $1< p\leq 2$,
\begin{align}
\norm{(1-\chi_{cr})\etag}_{B_{x_3}L_x^{p}}&\lesssim T^{\frac{m-1}{2}}r^{-m+\frac{2-p}{p}}
\label{estetagfarBz}\\
\norm{(1-\chi_{cr})\etacone}_{B_{x_3}L_x^{p}}&\lesssim T^{\frac{m-1}{2}}r^{-m+\frac{2-p}{p}}\nw
\label{estetaronefarBz}\\
\norm{(1-\chi_{cr})Q_{\Phi}^{-1}(\chitil_{8r}\wctwo)}_{B_{x_3}L_x^{p}}&\lesssim M^{-1}T^{\frac{m-1}{2}}r^{-m+\frac{2-p}{p}}\nw \label{estetartwofarBz}\\
\norm{(1-\chi_{cr})\etag}_{L^2_{x}}&\lesssim T^{\frac{m-1}{2}}r^{-m} \label{estetagL2} \\	
 \norm{(1-\chi_{cr})\etacone}_{L^2_{x}}&\lesssim T^{\frac{m-1}{2}}r^{-m}\nw \label{estetaroneL2}
\end{align}
\end{enumerate}
\end{lem}

\bpf
\begin{enumerate}
    \item \eqref{estetagBz} follows from the explicit formula \eqref{defetagvg}. \eqref{estetar1Bz} follow from \eqref{normdefinition}, Hölder inequality \eqref{lem:HolderBzLp} and weight estimate \eqref{estweightfull}. \eqref{estetar2Bz} is similar, using that 
     $$\norm{\etactwo}_{B_{x_3}L^{p}_x}\leq \norm{\chi_{2r}}_{B_{x_3}L^{\infty}_x}\norm{Q_\Phi^{-1}(\chitil_{8r}\wctwo)}_{B_{x_3}L_x^{p}}.$$
    For \eqref{estetar1r2L2w} use $B_{x_3}\subset L^\infty \subset L^p$ and \eqref{normdefinition}.
    \item For \eqref{estetagfarBz}, \eqref{estetaronefarBz} and \eqref{estetartwofarBz} proceed as before, using now \eqref{estweightfar}. \eqref{estetagL2} and \eqref{estetaroneL2} are identical, using again $B_{x_3}\subset L^2$.
\end{enumerate}
\epf
\begin{lem}[Vorticity bounds in physical frame]\label{lem:boundsphysicalframe} For $0<c\leq 6$,\leavevmode 
\begin{align}
			\norm{\chitil_{\frac{r}{4}}\wg}_{L^1_y}&\lesssim 1  \label{estwg1} \\
                \norm{\chitil_{\frac{r}{4}}\wcone}_{L^1_y}+M^{-1}\norm{\chitil_{\frac{r}{4}}\wctwo}_{L^1_y}&\lesssim \nw  \label{estwrL1}\\
			t^\frac{1}{2}\norm{\<t^{-\frac{1}{2}}\bd\>^m \wcone}_{L^2_{y}}&\lesssim \nw \label{estwrone2w}\\
			t^\frac{1}{2}\norm{\<t^{-\frac{1}{2}}\bd\>^m \chitil_{8r}\wctwo}_{L^2_{y}}
			+t^\frac{1}{2}\norm{\<t^{-\frac{1}{2}}\bd\>^m (1-\chitil_{6r})\wctwo}_{L^3_{y}}&\lesssim M^{-1}\nw \label{estwrtwo2w}\\
            \norm{(1-\chitil_{\frac{r}{4}})\wg}_{L^2_{y}}&\lesssim T^{\frac{m-1}{2}}r^{-m} \label{estwg2}\\
			\norm{(1-\chitil_{\frac{r}{4}})\wcone}_{L^2_{y}}&\lesssim T^{\frac{m-1}{2}}r^{-m}\nw \label{estwrone2}\\
			\norm{(1-\chitil_{cr})\wctwo}_{L^2_{y}}&\lesssim M^{-1}T^{\frac{m-1}{2}}r^{-m}\nw \label{estwrtwo2} \\
   \norm{(1-\chitil_{cr})\wctwo}_{L^{4/3}_{y}}&\lesssim M^{-1} r^{\frac{1}{2}-m}T^{\frac{m-1}{2}} \nw\label{estwrtwo43}
	\end{align}
\end{lem}
\bpf
To prove \eqref{estwg1} and \eqref{estwrL1} apply $Q_\Phi^{-1}$. Then, bound the change of variable \eqref{eq:defPQ} with Lemma \ref{changevariablebounds} and use Lemma \ref{lem:boundsstraightenedframe}. 
\eqref{estwrone2w}, \eqref{estwg2} and \eqref{estwrone2} follow similarly. Also the first term in \eqref{estwrtwo2w}, while the second is consequence of the definition \eqref{normdefinition}.
For \eqref{estwrtwo2} and \eqref{estwrtwo43} the strategy is: we split $$(1-\chitil_{cr})\wctwo=(1-\chitil_{cr})\chitil_{6r}\wctwo+(1-\chitil_{cr})(1-\chitil_{6r})\wctwo.$$
To bound the first part, apply $Q_\Phi^{-1}$ and use \eqref{estetartwofarBz}. To bound the second part, use Hölder, \eqref{normdefinition}, \eqref{estweightfar} and $r<1$.
\epf

\begin{lem}[Velocity bounds in straightened frame]\label{lem:estvelocitystraightened}\leavevmode 
	\begin{enumerate}
		\item \underline{Full space bounds}	\begin{align}
			t^\frac{1}{4}\norm{\vg}_{B_{x_3}L^4_x}+	t^\frac{1}{4}\norm{\na\vg}_{B_{x_3}L^{4/3}_x}&\lesssim 1 \label{estvg4}\\
			t^\frac{1}{4}\norm{\vcone}_{B_{x_3}L^4_x}+	t^\frac{1}{4}\norm{\na\vcone}_{B_{x_3}L^{4/3}_x}&\lesssim \nw \label{estvrone4}\\
			t^\frac{1}{4}\norm{\vctwo}_{B_{x_3}L^4_x}+	t^\frac{1}{4}\norm{\na\vctwo}_{B_{x_3}L^{4/3}_x}&\lesssim M^{-1}(1+T^{\frac{m-1}{2}}r^{-m-\frac{3}{4}})\nw \label{estvrtwo4}
		\end{align}
	\item \underline{Difference of velocities bounds} 
\begin{align}
		t^\frac{1}{4}\norm{\vg-\chi_r\vbarg}_{B_{x_3}L^4_x}+	t^\frac{1}{4}\norm{\na(\vg-\chi_r\vbarg)}_{B_{x_3}L^{4/3}_x}\lesssim T^\frac{1}{4}r^{-\frac{1}{2}}+r\ln{r^{-1}}+T^{\frac{m}{2}-\frac{1}{4}}r^{-m+\frac{1}{2}}\label{estdifvg4}\\
			t^\frac{1}{4}\norm{\vcone-(-\Delta)^{-1}\na\times\etacone}_{B_{x_3}L^4_x}	\lesssim (T^\frac{1}{4}r^{-\frac{1}{2}}+r\ln{r^{-1}}+T^{\frac{m}{2}-\frac{1}{4}}r^{-m+\frac{1}{2}})\nw \quad \label{estdifvrone4}
\end{align}

	\item \underline{Bounds far from filament}
	\begin{align}
		\norm{(1-\chi_r)\vbarg}_{B_{x_3}L^4_x}+	\norm{\na(1-\chi_r)\vbarg}_{B_{x_3}L^{4/3}_x}&\lesssim r^{-\frac{1}{2}}\label{estvbar4}\\
		\norm{(1-\chi_r)\vg}_{B_{x_3}L^4_x}+	\norm{P_\Phi^{-1}(\chitil_{7r}(1-\chitil_r)\ug)}_{B_{x_3}L^{4}_x}&\lesssim r^{-\frac{1}{2}}+T^{\frac{m-1}{2}}r^{-m+\frac{1}{2}}\label{estvug4}\\
		\norm{(1-\chi_r)\vcone}_{B_{x_3}L^4_x}+	\norm{P_\Phi^{-1}(\chitil_{7r}(1-\chitil_r)\ucone)}_{B_{x_3}L^{4}_x}&\lesssim (r^{-\frac{1}{2}}+T^{\frac{m-1}{2}}r^{-m+\frac{1}{2}})\nw\label{estvurone4}\\
		\norm{(1-\chi_r)\vctwo}_{B_{x_3}L^4_x}+	\norm{P_\Phi^{-1}(\chitil_{7r}(1-\chitil_r)\uctwo)}_{B_{x_3}L^{4}_x}&\lesssim M^{-1}(r^{-\frac{1}{2}}+r^{-m-\frac{1}{4}}T^{\frac{m-1}{2}})\nw\label{estvurtwo4}
	\end{align}

	\end{enumerate}
\end{lem}
\bpf
\begin{enumerate}
    \item For \eqref{estvg4} and \eqref{estvrone4} recall the definition of $v^g, \vctwo$ and apply \eqref{estBSeta} together with \eqref{estetagBz}, \eqref{estetar1Bz}. Obtaining \eqref{estvrtwo4} is more involving.

    By Sobolev embedding, we have that $\norm{\hat{v}^{c_2}}_{L^4_{x_1,x_2}}\leq \norm{\na_{x_1,x_2}\hat{v}^{c_2}}_{L^{4/3}_{x_1,x_2}}$, so is enough to bound $\norm{\na v^{c_2}}_{B_{x_3}L^{4/3}_{x}}$. We split into two terms.

    $$\na\vctwo=\underbrace{\na(P_\Phi^{-1}(\chitil_r(-\Delta)^{-1}\na\times(\chitil_{4r}\wctwo)))}_I+\underbrace{\na(P_\Phi^{-1}(\chitil_r(-\Delta)^{-1}\na\times((1-\chitil_{4r})\wctwo)))}_{II}.$$
    To control $I$ we apply \eqref{estBSeta} to $Q_\Phi^{-1}(\chitil_{4r}\wctwo)$ together with \eqref{estetar2Bz}. To deal with $II$, we notice that $II$ is supported on $|\widetilde{x}|\leq 2r$ and apply Lemma \ref{lem:bztolp}
    $$\norm{II}_{B_{x_3}L^{4/3}_x}\lesssim \norm{II}^{1/4}_{L^{4/3}_{x}}\norm{\pa_{x_3}II}^{3/4}_{L^{4/3}_{x}}+\norm{II}_{L^{4/3}_{x}}.$$
We begin estimating the term without $\pa_{x_3}$. Changing variables to cancel $P_\Phi^{-1}$ and applying gradient for products, 
$$\norm{II}_{L^{4/3}_{x}}\lesssim 
\norm{\na \chitil_r(-\Delta)^{-1}\na\times((1-\chitil_{4R})\wctwo)}_{L^{4/3}_{y}}
+\norm{\chitil_r\na(-\Delta)^{-1}\na\times((1-\chitil_{4r})\wctwo)}_{L^{4/3}_y}.$$
Notice that the supports of $\na\chitil_r\cup \chitil_r$ and $(1-\chitil_{4r})\wctwo$ are separated by $2r$. Denoting $K$ the Biot-Savart kernel in $\RR^2\times\TT$, we can write
$$
\begin{aligned}
\norm{II}_{L^{4/3}_{x}}\lesssim &
\norm{\na \chitil_r(((1-\chitil_{r})K)*((1-\chitil_{4r})\wctwo)))}_{L^{4/3}_{y}}\\
&+\norm{\chitil_r(((1-\chitil_{r})\na K)*((1-\chitil_{4r})\wctwo))}_{L^{4/3}_y}\\
\lesssim & \norm{\na \chitil_r}_{L^{2}_{y}}\norm{\boldone_{r\leq |y|}K}_{L^{2,\infty}_{y}}\norm{(1-\chitil_{4r})\wctwo}_{L^{4/3}_{y}}\\
&+\norm{\chitil_r}_{L^{4/3}_{y}}\norm{\boldone_{r\leq |y|}\na K}_{L^{2}_{y}}\norm{(1-\chitil_{4r})\wctwo}_{L^{2}_{y}}.
\end{aligned}$$
By direct computation we obtain $\norm{\na \chitil_r}_{L^{2}_{y}}\lesssim 1$, $\norm{\chitil_r}_{L^{4/3}_{y}}\lesssim r^{\frac{3}{2}}$. Lemma \ref{lem:BiotSavartbounds} gives $\norm{\boldone_{r\leq |y|} K}_{L^{2,\infty}_{y}}\lesssim r^{-\frac{1}{2}}$, $\norm{\boldone_{r\leq |y|}\na K}_{L^{2}_{y}}\lesssim r^{-\frac{3}{2}}$. Using \eqref{estwrtwo2} and \eqref{estwrtwo43}, we obtain
$$\norm{II}_{L^{4/3}_{x}}\lesssim M^{-1}T^{\frac{m-1}{2}}r^{-m}\nw.$$

The estimate of $\norm{\pa_{x_3} II}_{L^{4/3}_{x}}$ is similar. We obtain in this case
$$\norm{\pa_{x_3}II}_{L^{4/3}_{x}}\lesssim M^{-1}T^{\frac{m-1}{2}}r^{-m-1}\nw.$$
We conclude 
 $$\norm{II}_{B_{x_3}L^{4/3}_x}\lesssim M^{-1}T^{\frac{m-1}{2}}r^{-m-\frac{3}{4}}\nw.$$

 \item We start with \eqref{estdifvg4}. As we have explained before, it is enough to bound the term with gradient. We would like to apply \ref{prop:estdiffBS}, but we need to localize the support of $\etag$ first. We can split 

 $$
 \begin{aligned}
\na(\vg-\chi_r\vbarg)=&\na (\chi_r(P_\Phi^{-1}(-\Delta)^{-1}\na\times Q_{\Phi}-(-\Delta)^{-1}\na\times)\chi_{2r}\etag) \\
&+\na\chi_r(-\Delta)^{-1}\na\times(1-\chi_{2r})\etag \\
&+\chi_r(-\Delta)^{-1}\na^2\times(1-\chi_{2r})\etag
\end{aligned}$$
The first term is bounded directly with \eqref{estdiffBS} and \eqref{estetagBz}. For the second term, we apply lemmas \ref{lem:HolderBzLp}, \ref{lem:BSBzLp} and estimate \eqref{estetagfarBz}. For the third one, lemmas \ref{lem:HolderBzLp}, \ref{lem:RieszBzLp} and estimate \eqref{estetagfarBz}. 
For \eqref{estdifvrone4} we follow similar steps. We split in two terms
$$\vcone-(-\Delta)^{-1}\na\times\etacone=(\vcone-(-\Delta)^{-1}\na\times\chi_{2r}\etacone)+(-\Delta)^{-1}\na\times(1-\chi_{2r})\etacone.$$
The first term can be bounded using \eqref{estdiffBS} and the definition of $\nw$. The second term uses lemmas \ref{lem:HolderBzLp}, \ref{lem:BSBzLp} and estimate \eqref{estetaronefarBz}.
\item 
The estimate \eqref{estvbar4} follows from the definition \eqref{defetagvg}. For the other estimates, from \eqref{eq:relationuv} we see that the second terms are greater than the first ones. To bound \eqref{estvurone4} we want to apply Proposition \ref{prop:estBSeta}. To this end, we split in two terms 
$$
\begin{aligned}
P_\Phi^{-1}(\chitil_{7r}(1-\chitil_r)\ucone)=&P_\Phi^{-1}(\chitil_{7r}(1-\chitil_r)(-\Delta)^{-1}\na\times Q_\Phi((\chi_{2r}-\chi_{\frac{r}{4}})\etacone)\\
&+P_\Phi^{-1}(\chitil_{7r}(1-\chitil_r)(-\Delta)^{-1}\na\times Q_\Phi(\chi_{\frac{r}{4}}\etacone).
\end{aligned}
$$
The first term can be controlled using estimates \eqref{estBSeta} and \eqref{estetaronefarBz}. The second term can be controlled using the estimates \eqref{estBSawayeta} and \eqref{estetar1Bz}. The bound of \eqref{estvug4} is analogue, now using \eqref{estetagfarBz} and \eqref{estetagBz}.
For \eqref{estvurtwo4} we need extra care with the supports. Let $\rho(x)$ be smooth, radial, nonnegative, equal to $1$ on $\{|\widetilde{x}|\leq \frac{29}{2}r\}$ and null on $\{|\widetilde{x}|\geq \frac{31}{2}r\}$. Denote $\rhotil\circ\Phi \coloneqq \rho$. Note that $\chi_{7r}\lesssim \rho\lesssim \chi_{8r}$. We split 
\begin{equation}\label{eq:auxiliarstep}
\begin{aligned}
P_\Phi^{-1}(\chitil_{7r}(1-\chitil_r)\uctwo)=&P_\Phi^{-1}(\chitil_{7r}(1-\chitil_r)(-\Delta)^{-1}\na\times (\rhotil\wctwo)\\
&+P_\Phi^{-1}(\chitil_{7r}(1-\chitil_r)(-\Delta)^{-1}\na\times((1-\rhotil) \wctwo).
\end{aligned}
\end{equation}
We treat the first term as we did with \eqref{estvurone4}, i.e., we split 
$$
\begin{aligned}  
P_\Phi^{-1}(\chitil_{7r}(1-\chitil_r)(-\Delta)^{-1}\na\times (\rhotil\wctwo)=&P_\Phi^{-1}(\chitil_{7r}(1-\chitil_r)(-\Delta)^{-1}\na\times ((\rhotil-\chitil_{\frac{r}{4}})\wctwo) \\
&+P_\Phi^{-1}(\chitil_{7r}(1-\chitil_r)(-\Delta)^{-1}\na\times (\chitil_{\frac{r}{4}}\wctwo).
\end{aligned}
$$
We can apply Proposition \ref{prop:estBSeta} and estimate \eqref{estetartwofarBz} to the first term, getting 
$$
\begin{aligned}
\norm{P_\Phi^{-1}(\chitil_{7r}(1-\chitil_r)(-\Delta)^{-1}\na\times ((\rhotil-\chitil_{\frac{r}{4}})\wctwo)}_{B_{x_3}L^4_x}\lesssim \norm{Q_\Phi^{-1}((\rho-\chitil_{\frac{r}{4}})\wctwo)}_{B_{x_3}L^{4/3}_x} \\
\lesssim \norm{(1-\chi_{\frac{r}{4}})Q_\Phi^{-1}(\chitil_{8r}\wctwo)}_{B_{x_3}L^{4/3}_x}\\
\lesssim M^{-1}T^{\frac{m-1}{2}}r^{-m+\frac{1}{2}}\nw,
\end{aligned}
$$
while the second term can be controlled using \eqref{estBSawayeta} and \eqref{estetar2Bz}. To finish \eqref{estvurtwo4} we have to bound the second term of \eqref{eq:auxiliarstep}. This is analogue to the estimate of the term $II$ in the proof of \eqref{estvrtwo4}, using the fact that the supports of $\chitil_{7r}$ and $(1-\rhotil)$ are separated by $r/2$. This gives the estimate $M^{-1}r^{-m-\frac{1}{4}}T^{\frac{m-1}{2}}\nw$.
\end{enumerate}
\epf

\begin{lem}[Velocity bounds in physical frame]\label{lem:severalestimatesu}\leavevmode 
	\begin{align}
		r^\frac{1}{4}\norm{(1-\chitil_r)\ug}_{L^{12}_y}+	\norm{(1-\chitil_r)\ug}_{L^{6}_y}+\norm{\na(1-\chitil_r)\ug}_{L^{2}_y}&\lesssim r^{-\frac{3}{2}}+T^{\frac{m-1}{2}}r^{-m-\frac{1}{3}}\label{estugseveral}\\
		r^\frac{1}{4}\norm{(1-\chitil_r)\ucone}_{L^{12}_y}+	\norm{(1-\chitil_r)\ucone}_{L^{6}_y}+\norm{\na(1-\chitil_r)\ucone}_{L^{2}_y}&\lesssim (r^{-\frac{3}{2}}+T^{\frac{m-1}{2}}r^{-m-\frac{1}{3}})\nw \label{esturoneseveral}\\
	r^\frac{1}{4}\norm{(1-\chitil_r)\uctwo}_{L^{12}_y}+	\norm{(1-\chitil_r)\uctwo}_{L^{6}_y}+\norm{\na(1-\chitil_r)\uctwo}_{L^{2}_y}&\lesssim M^{-1}(r^{-\frac{3}{2}}+T^{\frac{m-1}{2}}r^{-m-\frac{1}{3}})\nw \label{esturtwoseveral}
	\end{align}
\end{lem}
\bpf
Estimates \eqref{estugseveral} and \eqref{esturoneseveral} are similar, we focus on the second one. We begin with the bound of the gradient. Again, we split 
$$(1-\chitil_r)\ucone=\underbrace{(1-\chitil_r)(-\Delta)^{-1}\na\times Q_\Phi((\chi_{2r}-\chi_{\frac{r}{4}})\eta^{c_1})}_I
+\underbrace{(1-\chitil_r)(-\Delta)^{-1}\na\times Q_\Phi(\chi_{\frac{r}{4}}\eta^{c_1})}_{II}.$$
We can bound 
$$
\begin{aligned}
\norm{\na I}_{L^2_y}\lesssim& \norm{\na(1-\chitil_r)}_{L^\infty_y}\norm{\chitil_{2r}}_{L^{3}_y}\norm{\boldone_{|y|\leq 8r}K}_{L^{\frac{3}{2},\infty}_y}\norm{Q_\Phi((\chi_{2r}-\chi_{\frac{r}{4}})\eta^{c_1})}_{L^2_y} \\
&+\norm{\na(-\Delta)^{-1}\na\times Q_\Phi((\chi_{2r}-\chi_{\frac{r}{4}})\eta^{c_1})}_{L^2_y}.
\end{aligned}
$$
The first term comes from applying gradient to $(1-\chitil_r)$, introducing an extra $\chitil_{2r}$ because it is identically $1$ on the support of the term, and using Hölder and Young inequality. We can localize the support of the Biot-Savart kernel because the support of other terms are also localized. Now, we bound these terms using the boundedness of the Riesz transform, Lemma \ref{changevariablebounds} and \eqref{estetaroneL2} to obtain
$$\norm{\na I}_{L^2_y}\lesssim (r^{-1}r^{\frac{2}{3}}T^{\frac{m-1}{2}}r^{-m}+T^{\frac{m-1}{2}}r^{-m})\nw\lesssim T^{\frac{m-1}{2}}r^{-m-\frac{1}{3}}\nw.$$
Using similar ideas, we can bound $II$ by 
$$
\begin{aligned}
\norm{\na II}_{L^2_y}\lesssim \norm{\na(1-\chitil_r)}_{L^\infty_y}\norm{\boldone_{\frac{r}{2}\leq |y|}K}_{L^{2,\infty}_y}\norm{Q_\Phi(\chi_{\frac{r}{4}}\eta^{c_1})}_{L^1_y}\lesssim r^{-1}r^{-\frac{1}{2}}\nw,
\end{aligned}
$$
where we have used \eqref{estwrL1}.
The estimate of the gradient term in \eqref{estugseveral} is analogue, using \eqref{estetagL2} and \eqref{estwg1}. Now we prove the $L^{12}$ estimate in \eqref{esturoneseveral}. We split in three terms
$$
\begin{aligned}
(1-\chitil_r)\ucone=&(1-\chitil_r)\chitil_{6r}(-\Delta)^{-1}\na\times Q_\Phi((\chi_{2r}-\chi_{\frac{r}{4}})\eta^{c_1}) \\
&+(1-\chitil_{6r})(-\Delta)^{-1}\na\times Q_\Phi((\chi_{2r}-\chi_{\frac{r}{4}})\eta^{c_1})
+(1-\chitil_r)(-\Delta)^{-1}\na\times Q_\Phi(\chi_{\frac{r}{4}}\eta^{c_1}).
\end{aligned}
$$
The third term can be bounded by
$$\norm{(1-\chitil_r)(-\Delta)^{-1}\na\times Q_\Phi(\chi_{\frac{r}{4}}\eta^{c_1})}_{L^{12}_y}\lesssim \norm{\boldone_{\frac{r}{2}\leq |y|} K }_{L^{12}_y}\norm{Q_\Phi(\chi_{\frac{r}{4}}\eta^{c_1})}_{L^1_y}\lesssim r^{-\frac{7}{4}}\nw.$$
For the second term, Young inequality and change of variables followed by Hölder inequality gives 
$$
\begin{aligned}
\norm{(1-\chitil_{6r})(-\Delta)^{-1}\na\times Q_\Phi((\chi_{2r}-\chi_{\frac{r}{4}})\eta^{c_1})}_{L^{12}_y}\lesssim &\norm{\boldone_{2r\leq |y| } K }_{L^{2,\infty}_y}\norm{Q_\Phi((\chi_{2r}-\chi_{\frac{r}{4}})\eta^{c_1})}_{L^\frac{12}{7}_y} \\
\lesssim & r^{-\frac{1}{2}}\norm{(\chi_{2r}-\chi_{\frac{r}{4}})\eta^{c_1}}_{L^2_{x}}\norm{\chi_{4r}}_{L^{12}_{x}} \\
\lesssim & r^{-\frac{1}{2}} T^{\frac{m-1}{2}}r^{-m}\nw r^{\frac{1}{6}}
\end{aligned}
$$
Finally, we bound the first term using Lemma \ref{lem:boundfourierK} by
$$
\begin{aligned}
\norm{\chitil_{6r}(-\Delta)^{-1}\na\times Q_\Phi((\chi_{2r}-\chi_{\frac{r}{4}})\eta^{c_1})}_{L^{12}_y}&\lesssim 
\norm{\FF Q_\Phi^{-1}\chitil_{6r}(-\Delta)^{-1}\na\times Q_\Phi((\chi_{2r}-\chi_{\frac{r}{4}})\eta^{c_1})}_{l^{\frac{12}{11}}_\zeta L^{12}_{x}} \\
&\lesssim \norm{\widetilde{K}}_{l^{\frac{12}{5}}_\zeta L^{2}_{x}}\norm{Q_\Phi((\chi_{2r}-\chi_{\frac{r}{4}})\eta^{c_1})}_{L^2_{x_3}L^\frac{12}{7}_x} \\
&\lesssim \norm{\chi_{4r}}_{L^{12}_{x}}\norm{(1-\chi_{\frac{r}{4}})\eta^{c_1})}_{L^2_{x}}\lesssim T^{\frac{m-1}{2}}r^{-m+\frac{1}{6}}\nw.
\end{aligned}$$
Note that we have obtained a better exponent in $r$ than the stated $T^{\frac{m-1}{2}}r^{-m-\frac{1}{3}}$, but that is enough for our purposes. The bound of $L^6$ norm follows in the same way than the $L^{12}$ norm. Also, \eqref{estugseveral} is analogue to \eqref{esturoneseveral}. The estimates for \eqref{esturtwoseveral} are harder because the support of $\wctwo$ is the full space. In compensation, we can play with having two different norms in \eqref{normdefinition}. For the estimate of $\na\uctwo$, we split 

$$(1-\chitil_r)\ucone=\underbrace{(1-\chitil_r)(-\Delta)^{-1}\na\times (1-\chi_{\frac{r}{4}})\wctwo}_I
+\underbrace{(1-\chitil_r)(-\Delta)^{-1}\na\times \chi_{\frac{r}{4}}\wctwo}_{II}.$$
Using estimates \eqref{estwrtwo2} and \eqref{estwrtwo43} and the boundedness of Riesz transform, 
$$
\begin{aligned}
\norm{\na I}_{L^2_y}\lesssim & \norm{\na(1-\chitil_r)}_{L^\infty_y}\norm{\chitil_{2r}}_{L^{4}_y}\norm{\boldone_{1\leq |y|}K}_{L^{2,\infty}_y}\norm{(1-\chi_{\frac{r}{4}})\wctwo}_{L^{\frac{4}{3}}_y} \\
&+\norm{\na(1-\chitil_r)}_{L^\infty_y}\norm{\chitil_{2r}}_{L^{3}_y}\norm{\boldone_{|y|\leq 1}K}_{L^{\frac{3}{2},\infty}_y}\norm{(1-\chi_{\frac{r}{4}})\wctwo}_{L^{2}_y} \\
&+\norm{\na(-\Delta)^{-1}\na\times(1-\chi_{\frac{r}{4}})\wctwo}_{L^2_y} \\
\lesssim & (r^{-1}R^{\frac{2}{3}}T^{\frac{m-1}{2}}r^{-m}+T^{\frac{m-1}{2}}r^{-m})M^{-1}\nw\lesssim T^{\frac{m-1}{2}}r^{-m-\frac{1}{3}}M^{-1}\nw.
\end{aligned}
$$
Using now \eqref{estwrL1},

$$
\begin{aligned}
\norm{\na II}_{L^2_y} &\lesssim \norm{\na(1-\chitil_r)}_{L^\infty_y}\norm{\boldone_{\frac{r}{2}\leq |y|\leq \frac{5r}{2}}K}_{L^2_y}\norm{\chitil_{\frac{r}{4}}\wctwo}_{L^1_y} 
+\norm{\boldone_{\frac{r}{2}\leq |y|}\na K }_{L^2_y}\norm{\chitil_{\frac{r}{4}}\wctwo}_{L^1_y} \\
&\lesssim M^{-1}r^{-\frac{3}{2}}\nw.
\end{aligned}
$$
For the $L^{12}$ bound of $\uctwo$, we split

$$(1-\chitil_r)\uctwo=(1-\chitil_r)\na\times(-\Delta)^{-1}\chitil_{6r}\wctwo+(1-\chitil_r)\na\times(-\Delta)^{-1}(1-\chitil_{6r})\wctwo.$$
The first term can be estimated similarly to $\norm{\ucone}_{L^{12}_y}$. For the second term, is enough to bound 
$$
\begin{aligned}
\norm{(1-\chitil_r)\na\times(-\Delta)^{-1}(1-\chitil_{6r})\wctwo}_{L^{12}_y}\lesssim &\norm{\boldone_{|y|\leq 1} K}_{L^{\frac{3}{2},\infty}_y}\norm{(1-\chitil_{6r})\wctwo}_{L^\frac{12}{5}_y} \\
&\hspace{-1cm}+\norm{\boldone_{1\leq |y|} K}_{L^{2,\infty}_y}\norm{(1-\chitil_{6r})\wctwo}_{L^\frac{12}{7}_y}\lesssim T^{\frac{m-1}{2}}r^{-m+\frac{1}{6}}\norm{\wctwo}_F,
\end{aligned}
$$
where we have used \eqref{estweightfar} and Hölder inequality. Again, the $L^6$ bound can be proved in a similar way than the $L^{12}$ bound.
\epf
\subsection{Estimates on $\etacone$}
In this subsection we will bound \eqref{eq:a1}. The following Lemma is a direct conclusion of \cite[Lemmas~7.9,7.10]{bedrossiangermainharropgriffiths23}. The proof is rather long and technical, so we will omit it. 

\begin{lem}\label{lemma:operatorSestimate}
	For any $0<s<t\leq T$ and $0\leq \beta<1/2$,
	\begin{equation*}
	\sqrt{t}\norm{\left<\frac{\widetilde{x}}{\sqrt{t}}\right>^m\<\sqrt{t}\na\>^\beta\int_0^t S(t,s)B[v,\eta]}_{B_{x_3}L^2_x}\lesssim t^\frac{1}{4}\norm{v}_{B_{x_3}L^4_x}\norm{\eta}_{N^0}.
	\end{equation*}
\end{lem}
Applying this lemma to equation \eqref{eq:a1}, we obtain

\begin{prop}\label{prop:boundsa1N} Assume $\eps<1, M>1, T \leq r^{22}$. If $\w\in B_{\eps,T,r,M}$, then 
$$\norm{a^{1}(\w)}_{N^\beta}\lesssim r\ln{r^{-1}}+\left(r\ln{r^{-1}}+M^{-1}\right)\eps+\eps^2.$$
If $\w, \w'$ have the same initial data, then 
$$\norm{a^{1}(\w)-a^{1}(\w')}_{N^\beta}\lesssim \left(r\ln{r^{-1}}+M^{-1}+\eps\right)\eps.$$
\end{prop}
 \bpf
The assumptions help to simplify the final bound. From \eqref{eq:a1} and Lemma \ref{lemma:operatorSestimate}, we have 
$$
\begin{aligned}
\norm{a^1}_{N^\beta}
\lesssim& t^\frac{1}{4}\norm{v-\vbarg}_{B_{x_3}L^4_x}\norm{\etacone}_{N^0}
+t^\frac{1}{4}\norm{v}_{B_{x_3}L^4_x}\norm{\etactwo}_{N^0}
+t^\frac{1}{4}\norm{\vg-\chi_r\vbarg}_{B_{x_3}L^4_x}\norm{\etag}_{N^0} \\
&+t^\frac{1}{4}\norm{v-\vg-(-\Delta)^{-1}\nabla\times\etacone}_{B_{x_3}L^4_x}\norm{\etag}_{N^0}.
\end{aligned}
$$
To control the first term, split 
$$v-\vbarg=(1-\chi_r)\vbarg+(\vg-\chi_r\vbarg)+\vcone+\vctwo$$
and use \eqref{estvbar4}, \eqref{estdifvg4}, \eqref{estvrone4}, \eqref{estvrtwo4} and \eqref{normdefinition}. For the second term, split $v$ and use \eqref{estvg4}, \eqref{estvrone4}, \eqref{estvrtwo4} and \eqref{normdefinition}. In the third term we apply \eqref{estdifvg4} and \eqref{estetagBz}. For the last term, we use \eqref{estdifvrone4},\eqref{estvrtwo4} and \eqref{estetagBz}. 

To bound the difference of solutions, notice that the terms were only $f^g$ variables appears cancel, since they don't depend on the solution. Then use that $a^1$ is linear and subtract similar terms. For example, 
$$B[v_1,\etactwo_1]-B[v_2,\etactwo_2]=B[v_1-v_2,\etactwo_1]-B[v_2,\etactwo_2-\etactwo_1].$$
Then we can proceed exactly as we did with the bound of $\norm{a^1}_{N^\beta}$.
 \epf
 
 \subsection{Estimates on $\wctwo$}
 In this subsection we bound \eqref{eq:a2}. We use two different norms for $\wctwo$, so we need control on both. 
\begin{prop}\label{prop:boundsa2N}Assume $\eps<1, M>1, T \leq r^{22}$. If $\w\in B_{\eps,T,r,M}$, then 
	$$M\norm{Q_\Phi^{-1}[\chitil_{8r}a^2(\w)]}_{N^0}\lesssim Mr+(Mr+r^\frac{1}{2})\eps+M\eps^2.$$
 If $\w, \w'$ have the same initial data, then 
$$M\norm{Q_\Phi^{-1}[\chitil_{8r}(a^2(\w)-a^2(\w'))]}_{N^0}\lesssim (Mr+r^\frac{1}{2}+M\eps)\eps.$$
\end{prop}
\bpf 
Checking \eqref{eq:a2}, the different contributions to $a^2$ come from $B[(1-\chitil_r)u,\w]+\EE^{g}+\EE^{c}$. 

From their definition \eqref{eq:evwg}\eqref{eqevwrone} and Lemma \ref{lemma:operatorstransformations}, notice that both $Q_\Phi^{-1}(\EE^g)$ and $Q_\Phi^{-1}(\EE^c)$ are supported on $\{|\widetilde{x}|\leq 4r\}$ and can be written as 
\begin{equation}\label{eq:expressionEE}
Q_\Phi^{-1}(\EE^*)=\chi_{2r}\Delta(\eta^*)-\Delta_\Phi(\chi_{2r}\eta^*)=\na^\gamma_{x,z}(C_\gamma\eta^*),
\end{equation}  
where all matrix $C_\gamma$ are smooth and satisfy, by Lemma \ref{changevariablebounds},
$$|C_2|\lesssim r, \quad |C_1|\lesssim r^{-1}, \quad |C_0|\lesssim r^{-2}.$$
By Lemma \ref{lem:boundduhamelheat}, \eqref{normdefinition} and \eqref{estetagBz}, the term with $\EE^g+\EE^c$ is bounded by $(r+T^\frac{1}{2}r^{-1}+Tr^{-2})(1+\nw)$. 
Now, we split $B[(1-\chitil_r)u,\w]=B[(1-\chitil_r)u^g,\w]+B[(1-\chitil_r)u^r,\w]$. To control the contribution of $B[(1-\chitil_r)u^g,\wg+\wcone+\chitil_{6r}\wctwo ]$ we use Lemmas \ref{lem:boundduhamelheat} with $f_2\equiv 0$ and \ref{lem:HolderBzLp} to get
$$
\begin{aligned}
& \norm{Q_\Phi^{-1}\left(\chitil_{8r}\int_0^te^{(t-s)\Delta}B[(1-\chitil_r)u^g,\wg+\wcone+\chitil_{6r}\wctwo]ds\right)}_{N^0}\\ 
&\hspace{3cm} \lesssim  t^\frac{3}{4} \norm{P_\Phi^{-1}(\chitil_{7r}(1-\chitil_r)u^g)}_{B_{x_3}L^4_x}
\norm{\< t^{-\frac{1}{2}}x\>^m(\etag+\etacone+Q_\Phi^{-1}(\chitil_{6r}\wctwo))}_{B_{x_3}L^2_x} \\
&\hspace{3cm} \lesssim T^{\frac{1}{4}}r^{-\frac{1}{2}}(1+\nw+M^{-1}\nw),
\end{aligned}
$$
where we used for the last inequality \eqref{estvug4}, \eqref{estetagBz} and \eqref{normdefinition}. For the term with $B[(1-\chitil_{r})u^g,(1-\chitil_{6r})\wctwo]$ we apply Lemma \ref{lem:boundduhamelheat} with $f_1\equiv 0$ followed by Hölder inequality, \eqref{estugseveral} and \eqref{normdefinition} to obtain
$$\norm{Q_\Phi^{-1}\left(\chitil_{8r}\int_0^te^{(t-s)\Delta}B[(1-\chitil_r)u^g,(1-\chitil_{6r})\wctwo]ds\right)}_{N^0}\lesssim T^\frac{1}{4}r^{-\frac{3}{2}}M^{-1}\nw.$$

The bound for $B[(1-\chitil_r)u^r,\w]$ is obtained in the same way, using \eqref{estvurone4} and \eqref{estvurtwo4} instead of \eqref{estvug4}, and \eqref{esturoneseveral}, \eqref{esturtwoseveral} instead of \eqref{estugseveral}.
The estimate for $a^2(\w)-a^2(\w')$ follows from a similar procedure.
\epf

\begin{prop}\label{prop:boundsa2F}
Assume $\eps<1, M>1, T \leq r^{22}$. If $\w\in B_{\eps,T,r,M}$, then  
	$$M\norm{a^2}_F\lesssim Mr+Mr\eps+M\eps^2.$$
 If $\w, \w'$ have the same initial data, then 
 $$M\norm{a^2(\w)-a^2(\w')}_F\lesssim (Mr+M\eps)\eps.$$
\end{prop}
We start with two auxiliary lemmas

\begin{lem}\label{lem:bounddistanceheat}
With the notation of Lemma \ref{lem:Heatproductspace},
$$\boldone_{r\leq |y|}H_\RRT(y)\lesssim e^{-\frac{r^2}{16t}}e^{-\frac{|y|^2}{16t}}t^{-\frac{3}{2}}e^{-\frac{|\widetilde{y}|^2}{8t}}\sum_{k\in \ZZ}e^{-\frac{|y_3+2k\pi|^2}{8t}},$$
i.e., we can subtract a factor $e^{-\frac{r^2+|y|^2}{16t}}$ if we change the denominator of the exponentials defining the heat kernel.
\end{lem}
\bpf
Use triangular inequality and notice that $e^{-\frac{|y_3+2k\pi|^2}{4t}}\leq e^{-\frac{|y_3|^2}{8t}}e^{-\frac{|y_3+2k\pi|^2}{8t}}$ for all $k\in \ZZ$. Do the same with $e^{-\frac{|\widetilde{y}|^2}{4t}}$, and use $r\leq |y|$.
\epf

\begin{lem}\label{lem:weightheat}
For any $y,y'\in \RRT$ and $0<s\leq t$,
$$\left< \frac{\bd(y)}{\sqrt{t}}\right>^me^{-\frac{|y-y'|^2}{16t}}\lesssim \left< \frac{\bd(y')}{\sqrt{s}}\right>^m$$
\end{lem}
\bpf
Consider the cases $|y|\geq 2|y'|$ and $|y|\leq 2|y'|$. In the first case, the LHS is bounded because $|\bd(y)|\to\infty\iff |y|\to\infty$ and the exponential has faster decay. In the second case, $\bd(y)\approx \bd(y')$.
\epf

Now we proceed with the proof of Proposition \ref{prop:boundsa2F}.
\bpf
We begin with the bound for $a^2$.

For the term with $\EE^c$, use the formula \eqref{eq:expressionEE}. We can pass derivatives through $Q_\Phi$ by adding lower-order terms, due to equation \eqref{eq:defPQ} and Lemma \ref{changevariablebounds}. Then pass the derivatives to the heat kernel. Then, use the separation of supports of $C_\gamma, (1-\chitil_{6r})$ to apply Lemma \ref{lem:bounddistanceheat}, and \ref{lem:weightheat} to pass the weight to the vorticity term. Then, using the Young inequality, we obtain

$$\norm{(1-\chitil_{6r})\< t^{-\frac{1}{2}}\bd\>^m\int_0^te^{(t-s)\Delta}Q_\Phi(\na^\gamma_{x}C_\gamma\etacone)}_{L^3_y}\lesssim \int_0^t (t-s)^{-\frac{1+2|\gamma|}{4}}e^{-\frac{r^2}{16(t-s)}}s^{-\frac{1}{2}}\norm{C_\gamma}_{L^\infty}\nw ds.$$
Exponential can control any power of $(t-s)$, so at the end we get 

$$\norm{\int_0^te^{(t-s)\Delta}\EE^c}_F\lesssim T r^{-2} \nw.$$
The term with $\EE^g$ is similar. 

We now treat the term with $B[(1-\chitil_r)u^g,\wctwo]$. Begin by splitting in

$$\underbrace{\norm{\int_0^t e^{(t-s)\Delta}B[(1-\chitil_r)u^g, \chitil_{6r}\wctwo]}_F}_I+\underbrace{\norm{\int_0^t  e^{(t-s)\Delta}B[(1-\chitil_r)u^g,(1-\chitil_{6r})\wctwo]}_F}_{II}.$$
The strategy to bound both terms will be the same, but using different $L^p$ norms. In both cases we pass the divergence to the heat kernel and apply Young and Hölder inequalities to get
$$
\begin{aligned}
I &\leq t^\frac{1}{2}\int_0^t \norm{\< s^{-\frac{1}{2}}\bd\>^m\divergence e^{(t-s)\Delta}}_{L^{\frac{4}{3}}_y}\norm{(1-\chitil_{r})u^g}_{L^{12}_y}\norm{\< s^{-\frac{1}{2}}\bd\>^m\chitil_{6r}\wctwo}_{L^2_y}ds \\
&\lesssim t^\frac{1}{2}\int_0^t t^{-\frac{11}{8}}r^{-\frac{7}{4}}\nw (1-s/t)^{-\frac{7}{8}}(s/t)^{-\frac{1}{2}}ds\lesssim t^\frac{1}{8}r^{-\frac{7}{4}}\nw\leq r \nw.
\end{aligned}
$$
For $II$, we get 
$$
\begin{aligned}
II &\leq t^\frac{1}{2}\int_0^t \norm{\< s^{-\frac{1}{2}}\bd\>^m\divergence e^{(t-s)\Delta}}_{L^{\frac{6}{5}}_y}\norm{(1-\chitil_{r})u^g}_{L^{6}_y}\norm{\< s^{-\frac{1}{2}}\bd\>^m(1-\chitil_{6r})\wctwo}_{L^3_y}ds \\
&\lesssim t^\frac{1}{2}\int_0^t t^{-\frac{5}{4}}r^{-\frac{3}{2}}\nw (1-s/t)^{-\frac{3}{4}}(s/t)^{-\frac{1}{2}}ds\lesssim t^\frac{1}{4}r^{-\frac{3}{2}}\nw.
\end{aligned}
$$
The treat of the term with $B[(1-\chitil_r)u^g,\wcone]$ is analogue to $I$. For $B[(1-\chitil_r)u^g,\wg]$ is the same, except that $\nw$ does not appear. The remaining term, $B[(1-\chitil_r)u^r,\w]$, is treated identically using the corresponding inequalities from Lemma \ref{lem:severalestimatesu}. 
The bound of $a^2(\w)-a^2(\w')$ is similar. 
\epf

\subsection{Short time existence}

\begin{theo}\label{theo:fixedpoint}
There exists $T, M, r, \eps$ such that $\QQ: B_{\eps,T,r,M}\to B_{\eps,T,r,M}$ and for $\w_1, \w_2 \in B_{\eps,T,r,M}$ there holds the contraction property $||\QQ(\w_1)-\QQ(\w_2)||_{X_t}\leq \frac{1}{2}||\w_1-\w_2||_{X_t}$. It follows that there exists a unique fixed point $\QQ(\w)=\w\in  B_{\eps,T,r,M}$. 
\end{theo}
\bpf  
The Theorem is a consequence of Propositions \ref{prop:boundsa1N}, \ref{prop:boundsa2N} and \ref{prop:boundsa2F}. Assume that $M>1, T<r^{22}, \eps<1$.  If $\w\in B_{\eps,T,r,M}$, then 
$$\norm{\QQ(\w)}_{X_T}\lesssim  r\ln{r^{-1}}+Mr+\eps(r\ln{r^{-1}}+M^{-1}+Mr+r^{1/2})+\eps^2M.$$
To get a bound smaller than $\eps$, take $M\gg  1, \eps\ll M^{-1}, r\ll \eps.$ This choice also makes $\QQ$ contractive, which implies uniqueness of the fixed point.
\epf

\begin{theo}\label{prop:fixedpointsolvesNS}
The solution constructed in Theorem \ref{theo:fixedpoint} is a solution of the 3D Navier-Stokes equations in the sense of Definition \ref{def:mildsolution}. Therefore, Theorem \ref{theo:localexistence} holds. 
\end{theo}
\bpf 
The fact that for positive times $\w\in L^1\cap L^2$ follows from the vorticity decomposition and the norms used \eqref{normdefinition}. The integrability for the velocity $u$ can be deduced from the Biot-Savart law.
The rest of the proof is analogue to \cite[Proposition 7.18]{bedrossiangermainharropgriffiths23}, so we skip it. 
\epf

\begin{coro}[Uniqueness]\label{coro:uniquenesssmalltime}
Let $\bar{\w}\in C^2((0,T^*))\times \RR^3$ be another mild solution such that for sufficiently small $r>0$ we have
\begin{equation}\label{eq:uniquenesscondition}
\lim_{T\to 0^+} \norm{Q_\Phi^{-1}(\chitil_{8r}(\bar{\w}-\wg))}_{N^0}+\norm{\bar{\w}-\wg}_F=0.
\end{equation}
Then $\bar{\w}=\w$, where $\w$ is the solution constructed in Theorem \ref{theo:fixedpoint}.
\end{coro}
\bpf
In the proof of Theorem \ref{theo:fixedpoint}, the last constants we choose are $r$ and $T$. We keep the choice of $M$ and $\eps$.
For any $r>0$ small enough, we can decompose $\bar{\w}$ in $\bar{\w}=\bar{\w^g}+\bar{\w}^{c_1}+\bar{\w}^{c_2}$, where $Q_\Phi^{-1}(\chitil_{2r}\bar{\w}^{c_1})$ is forced to satisfy equation \eqref{eqevetarone} and $\bar{\w}^{c_2}$ is the remainder. Note that this decomposition only depends on $r$. Then, use \eqref{eq:uniquenesscondition} to take a small enough $T$ so that the fixed-point theorem can be applied on the ball \eqref{eq:uniquenesscondition}, where we have uniqueness. For more details, see \cite[Theorem 7.19]{bedrossiangermainharropgriffiths23}.
\epf

\section{Preservation of helical symmetry}\label{Section:preservationsymmetry}
In this section we check the symmetry of the solution obtained in Theorem \eqref{theo:fixedpoint}, provided the initial data have helical symmetry, which is defined in Definition \ref{def:helicalsymmetry}. Due to rotation invariance of the Navier-Stokes equations, any single filament with helical symmetry can be reduced to $\gamma(x_3)=(\rho\cos(x_3),\rho\sin(x_3),\kappa x_3).$ Once fixed, the parameters $\rho, \kappa$ do not play any role. To simplify the exposition, we will assume $\rho=\kappa=1$, so that the initial vorticity is supported on the curve 

$$\gamma(x_3)=(\cos(x_3),\sin(x_3),x_3).$$
Navier-Stokes equations preserve helical symmetry \cite{MahalovTitiLeibovich90}, but we cannot apply this property directly because the initial data is not smooth. We could try to smooth the initial data and show convergence, but this approach does not work. More precisely, if we use time-independent mollifiers and denote the smoothed initial data in the straightened variables by $\eta_0^\eps$, one can check that $\norm{\eta_0-\eta_0^\eps}_{N^0}\nrightarrow
0$. We now check the conservation of symmetry in our construction.

\subsection{Helicity of $\wg$}
Following subsection \ref{subsection:geometricalsetting}, we obtain that 
$$\vect(x_3)=\frac{1}{\sqrt{2}}(-\sin(x_3),\cos(x_3),1), \quad \vecn(x_3)=(-\cos(x_3),-\sin(x_3),0), \quad \vecb(x_3)=\frac{1}{\sqrt{2}}(\sin(x_3),-\cos(x_3),1),
$$
so the mapping between the straightened and curved frame is
\begin{equation}\label{eq:phihelicalcase}
\Phi(x)=(\cos(x_3),\sin(x_3),x_3)+x_1(-\cos(x_3),-\sin(x_3),0)
+\frac{x_2}{\sqrt{2}}(\sin(x_3),-\cos(x_3),1),
\end{equation}
with jacobian
\begin{equation}\label{eq:Jhelicalcase}
J=\left[\vecn \quad \vecb \quad \left( \frac{2-x_1}{\sqrt{2}}\vect+\frac{-x_2}{\sqrt{2}}\vecn+\frac{x_1}{\sqrt{2}}\vecb \right) \right],
\end{equation}
and determinant $\fancyD=\frac{2-x_1}{\sqrt{2}}$. We prove that $\wg$ is helical in three steps.

\begin{lem}
The image by $\Phi:\Sigma_r\to\Gam_r$ of vertical lines are helical lines. The preimage of helical lines are vertical lines.
\end{lem}
\bpf All we need is to focus on equation \eqref{eq:phihelicalcase}. If we fix $x_1, x_2$, we obtain vectors that rotate in the horizontal plane with angle $x_3$ and a vertical translation with jump $x_3$. For the second statement find the preimage of a point, and use that $\Phi$ is bijective and the first statement of the Lemma. 
\epf

\begin{coro}
By the previous lemma, the distance from $\Phi(x)$ to the initial vortex curve $\Gam$ is constant as a function of $x_3$.
\end{coro}
\bpf The set of points $\Phi(x)$ with $\widetilde{x}$ fixed is a helix with the same symmetry than $\Gam$. Let $\underline{x}_3$ be a value of $x_3$ that minimizes the distance between $\Phi(x)$ and $\Gam$. Applying the operator $S_\theta$ defined in \eqref{eq:helicalStheta} we can easily see that the distance does not depend on $x_3$.

\epf

\begin{prop}\label{prop:zsymmetryiffhelical}
If $\eta$ does not depend on $x_3$, $Q_\Phi(\eta)$ has helical symmetry. Conversely, if $\w$ has helical symmetry, $Q_\Phi^{-1}(\w)$ does not depend on $x_3$.
\end{prop}
\bpf For the first statement we must check that, restricted to a helix, $Q_\Phi(\eta)$ rotates with vertical axis and angle $x_3$. Due to the previous lemma, in this restriction $x_1,x_2$ are fixed. $Q_\Phi$ is defined in \eqref{eq:defPQ}, where we see that everything on the right-hand side is constant except $J$ \eqref{eq:Jhelicalcase}, which rotates as needed. The second statement is completely analogue since our operators are invertible. 
\epf

\begin{coro}
If $\Gamma$ is helical, then $\wg$ defined in \eqref{def:wgubarg} has helical symmetry.
\end{coro}
\bpf 
Recall that $\wg$ is defined from $\etag$, and we see in \eqref{defetagvg} that $\etag$ does not depend on $x_3$. Note that multiplying by $\chi_{2r}$ or $\chitil_{2r}$ does not affect symmetry. By Proposition \ref{prop:zsymmetryiffhelical}, we conclude that $\wg$ is helical. 
\epf

\subsection{Helicity of $\w$}
\begin{prop}
Let $\w$ be the solution of Theorem \ref{theo:fixedpoint}, with helical initial data. Then $\w(t)$ is helical for all $0\leq t\leq T$. Therefore, Proposition \ref{prop:helicalsymmetrypreserved} holds.
\end{prop} 
\bpf Note that the Biot-Savart operator and the Laplacian applied to functions that are $x_3$-independent keep the $x_3$-independence. These operators also preserve helical symmetry. By Proposition \ref{prop:zsymmetryiffhelical}, $x_3$-independence of $v^*$ or $\eta^*$ and the helical symmetry of $u^*$ or $\w^*$ are equivalent. Also, the operator $S$ used in \eqref{eq:a1} preserves $x_3$-independence when the variables involved are $x_3$-independent. By the same reasoning, the $\EE^*$ appearing in \eqref{eq:a2} are helical if the corresponding variables are $x_3$-independent in the straightened frame, and helical in the physical frame. 

Recall that $\w=\wg+Q_\Phi(\chi_{2r}\etacone)+\wctwo$. At time $0$, $\etacone$ and $\wctwo$ are $0$, so initially they have the desired symmetry. By the previous reasoning, if $f^{c_1}$ and $f^{c_2}$ in \eqref{eq:a1} and \eqref{eq:a2} have the desired symmetry, then $a^1,a^2$ will also have it. Since each iteration will preserve the symmetry, the fixed point solution $\w$ has the desired symmetry.
\epf

\section{Global in time existence for helical initial data}\label{section:globaltimeexistence}

By the previous section, we know that the solution to the Navier-Stokes equations obtained in theorem \ref{theo:fixedpoint} is helical if $\Gam$ is helical. Our purpose now is to use this fact to extend the solution local-in-time to a global-in-time smooth solution. We cannot apply the techniques of \cite{MahalovTitiLeibovich90} because the solution obtained in \ref{theo:localexistence} belongs only to $\LTu$. We apply results from local energy weak solutions instead. First, we show that smooth helical functions with compact support are dense in the space of helical functions belonging to $\mathring{E}_2$. Omitting the symmetry restrictions, this result is well known. 

\begin{lem}\label{lem:approxlocalenergy}
For any $u\in \mathring{E}_2$ helical and for any $\eps>0$, there exist divergence free $u^\eps\in C^\infty_c(\RRT)$ helical such that 
$$\norm{u-u^\eps}_{\LTu}<\eps.$$
\end{lem}
\bpf
Following \cite[Lemma 6.1]{kikuchi2007weak} and using cylinders instead of balls, we obtain a function $v^R$ supported in $C(0,2R)$ such that 
\begin{equation}\label{DivergenceCorrectorEstimate}
\begin{aligned}
    \divergence v^R&=u\cdot \na \chi_R, \\
    \norm{v^R}_{\LTu}&\lesssim \frac{\norm{u}_{\LTu}}{R} .
\end{aligned}
\end{equation}
In the construction, it is not clear whether $v^R$ is helical or not. However, using that $\divergence v^R$ is helical, we can construct a helical version of $v^R$. Take cylindrical coordinates $(r,\theta,z)$ and make the change $\theta\to \xi=\theta+z$. Then define 
$$\bar{v}^R=\frac{1}{2\pi}\int_0^{2\pi} v^R(r,\xi,z+s)ds.$$
By construction, $\bar{v}^R$ is helical and also satisfies estimates \eqref{DivergenceCorrectorEstimate}.
We conclude that $u^R=u\chi_R-\bar{v}^R$ is divergence free and 
$$\norm{u-u^R}_{\LTu}\leq \norm{u-u\chi_R}_{\LTu}+\norm{\bar{v}^R}_{\LTu}<\eps$$
if $R\gg 1$. We mollify $u^R$ to get $u^\eps$.
\epf

\begin{theo}\label{theo:globalweaklocalenergysolutions}
Let $u_0\in \mathring{E}_2$. Then the Navier-Stokes equations \eqref{NSvelocity} have at least one global in time solution in the sense of \ref{def:localenergyweaksolutions}. Furthermore, if $u_0$ is helical, at least one solution is helical.
\end{theo}
\bpf
Without the helical part, this is proven in \cite{lemarierieusset2002} and also in \cite{kikuchi2007weak} in the $\RR^3$ case. Following \cite{kikuchi2007weak}, the adaptation to the $\RRT$ case can be made taking into account the different decay at infinity (which will be similar to 2D) and the change of domain. To prove the helical part we need to obtain a smooth helical approximation to the initial data, which is done in Lemma \ref{lem:approxlocalenergy}. The conservation of the symmetry in the solution is guaranteed since the strategy is to solve Navier-Stokes or a regularized version which preserves the symmetry for regular data, and the limit of helical functions is helical. The complete proof would be rather long, so we omit it.
\epf

Now, we will perform energy estimates to show that the helical solution constructed in Theorem \ref{theo:globalweaklocalenergysolutions} is unique and has improved regularity. We will need to estimate terms involving the pressure, for which we will use the following proposition.
 
\begin{prop}\label{prop:pressuredecomposition}
For any $x_0\in \RRT$ and $(t,x)\in (0,T)\times C(x_0,3/2)$, the pressure $p$ associated to the solution $u$ of Theorem \ref{theo:globalweaklocalenergysolutions} can be written as
\begin{equation}\label{eq:pressuredecomposition}
\begin{aligned}
p_{x_0}(t,x)&=p_{x_0,1}(t,x)+p_{x_0,2}(t,x), \\
p_{x_0,1}(t,x)&=-\frac{1}{3}|u(t,x)|^2+\frac{1}{4\pi}\int_{C(x_0,2)}N_{i,j}(x-y):u(t,y)\otimes u(t,y)dy, \\
p_{x_0,2}(t,x)&=\frac{1}{4\pi}\int_{(\RRT)\setminus C(x_0,2)}(N_{i,j}(x-y)-N_{i,j}(x_0,y)):u(t,y)\otimes u(t,y)dy.
\end{aligned}
\end{equation}
In particular, for any $q>1$ we have $\norm{p_{x_0,1}}_{L^q(C(x_0,1))}\lesssim \norm{u}_{L^{2q}_{\uloc}}^2,$ $\norm{p_{x_0,2}}_{L^\infty(C(x_0,1))}\lesssim \norm{u}_{L^{2}_{\uloc}}^2$. 
\end{prop}
\bpf
The decomposition comes from equation \eqref{eq:pressuredecomposition}. The Newtonian potential $N$ behaves as the potential in three dimensions close to $0$, and as the potential in two dimensions far from $0$. See Lemma \ref{lem:BiotSavartbounds} or \cite[Proposition 3.1]{bronzinussenzveig2015} for a different proof.

To prove the bound for $p_{x_0,1}$, we use theory of singular integral operators. In order to get the estimate for $p_{x_0,2}$, we first note that in the domain of the integral $|N_{i,j}(x-y)-N_{i,j}(x_0,y)|\lesssim \frac{1}{|x_0-y|^3}$. Then, we can bound 
$$
\begin{aligned}
|p(t,x)_{x_0,2}|\lesssim& \sum_{i=1}^\infty \int_{C(x_0,2^{i+1})\setminus C(x_0,2^{i})}\frac{1}{|x_0-y|^3}|u(t,y)|^2dy  \\
\lesssim& \sum_{i=1}^\infty \frac{1}{2^{3i}} \int_{C(x_0,2^{i+1})}|u(t,y)|^2dy \lesssim 
\sum_{i=1}^\infty \frac{1}{2^{3i}} 2^{2i} \norm{u}_{\LTu}^2= \norm{u}_{\LTu}^2.
\end{aligned}
$$
\epf 

So far, we have shown that the helical symmetry is preserved but we have not exploited it. We do that in the following lemma. 
\begin{lem}\label{lem:helicalembedding}
Let $u$ be a smooth helical vector field, $x=(x_1,x_2,x_3)\in \RRT$. 
Then
 
\begin{equation}\label{eq:helicalinterpolation}
	\norm{u}_{L^4(C(x,1))}^4\lesssim \norm{u}_{L^2(C(x,2))}^2\norm{\nabla u}_{L^2(C(x,2))}^2+\norm{u}_{L^2(C(x,2))}^4.
\end{equation}
\end{lem}
\bpf 
The idea is to use the symmetry to reduce the estimate to a two dimensional estimate, then recover the three dimensional domain. We use $(\rho, \theta, z)$ as cylindrical coordinates, 
$\xi=\theta+z$, and $x=(\rho_x,\theta_x,z_x$). 
Observe that if $(\rho_1,\theta_1,z_1), (\rho_2,\theta_2,z_2)\in C(x,1)$, then $|\theta_1-\theta_2|\lesssim \frac{1}{1+\rho_x}$. This implies that for fixed $\rho, \xi$; if $(\rho,\xi,z_1), (\rho,\xi, z_2)\in C(x,1)$, then $|z_1-z_2|\lesssim \frac{1}{1+\rho_x}$. 

We will use cylindrical $(\rho,\theta,z)$ and helical $(\rho,\xi,z)$ coordinates. Because of helical symmetry, in helical coordinates $u$ does not depend on $z$, so we can write $u(\rho,\theta,z)=\underline{u}(\rho,\xi)$. From this expression, we deduce $\pa_\theta u(\rho,\theta,z)=\pa_\xi\underline{u}(\rho,\xi).$

We bound
$$
\begin{aligned}
\norm{u}_{L^4(C(x,1))}^4&=\int_{\{(\rho,\theta,z)\in C(x,1)\}}|u(\rho,\theta,z)|^4\rho d\rho d\theta dz \lesssim \int_{\{(\rho,\xi,z)\in C(x,1)\}}|\underline{u}(\rho,\xi)|^4\rho d\rho d\xi dz\\
&\lesssim \dfrac{1}{1+\rho_x}\int_{-\pi}^\pi\int_{\{|\rho-\rho_x|<1\}}|\underline{u}(\rho,\xi)|^4\rho d\rho d\xi,
\end{aligned}
$$
The last integral can be seen as a two-dimensional norm in polar coordinates. By Gagliardo-Nirenberg inequality, it follows that
$$
\begin{aligned}
\norm{u}_{L^4(C(x,1))}^4&\lesssim \dfrac{1}{1+\rho_x}\int_{-\pi}^\pi\int_{\{|\rho-\rho_x|<1\}}|\underline{u}
(\rho,\xi)|^2\rho d\rho d\xi\int_{-\pi}^\pi\int_{\{|\rho-\rho_x|<1\}}(|\pa_\rho \underline{u}
(\rho,\xi)|^2+\frac{1}{\rho^2}|\pa_\xi \underline{u}
(\rho,\xi)|^2)\rho d\rho d\xi\\
&+\dfrac{1}{1+\rho_x}\left(\int_{-\pi}^\pi\int_{\{|\rho-\rho_x|<1\}}|\underline{u}
(\rho,\xi)|^2\rho d\rho d\xi\right)^2\\
&\hspace{-0.4cm}\lesssim \int_{\{(\rho,\xi,z)\in C(x,2)\}}|\underline{u}
(\rho,\xi)|^2\rho d\rho d\xi dz \int_{\{(\rho,\xi,z)\in C(x,2)\}}(|\pa_\rho \underline{u}
(\rho,\xi)|^2+\frac{1}{\rho^2}|\pa_\xi \underline{u}
(\rho,\xi)|^2)\rho d\rho d\xi dz \\
&+\left(\int_{\{(\rho,\xi,z)\in C(x,2)\}}|\underline{u}
(\rho,\xi)|^2\rho d\rho d\xi dz\right)^2 \\
&\hspace{-0.4cm}\lesssim  \int_{\{(\rho,\theta,z)\in C(x,2)\}}|u(\rho,\theta,z)|^2\rho d\rho d\theta dz \int_{\{(\rho,\theta,z)\in C(x,2)\}}(|\pa_\rho u
(\rho,\theta,z)|^2+\frac{1}{\rho^2}|\pa_\theta u
(\rho,\theta,z)|^2)\rho d\rho d\theta dz \\
&+\left(\int_{\{(\rho,\theta,z)\in C(x,2)\}}|u(\rho,\theta,z)|^2\rho d\rho d\theta dz\right)^2=\norm{u}_{L^2(C(x,2))}^2\norm{\nabla u}_{L^2(C(x,2))}^2+\norm{u}_{L^2(C(x,2))}^4.
\end{aligned}$$
In the second step, we used that fixed $\rho,\xi$ such that $|\rho-\rho_x|<1, \xi\in \TT$, then $\mu\{z:(\rho,\xi, z)\in C(x,2)\}>\frac{1}{\rho_x+1}$, where $\mu$ is the one-dimensional Lebesgue measure. This is equivalent to saying that, in cylindrical coordinates, if we fix $\rho,z$ such that $|\rho-\rho_x|<1$, then $\mu\{\theta:(\rho,\theta, z)\in C(x,2)\}>\frac{1}{\rho_x+1}$, which is more intuitive. Thus, we can recover the integral in $z$ by adding a factor $\rho_x+1$.
\epf
\begin{remark}
The estimate \eqref{eq:helicalinterpolation} is known in cylinders centered at $0$ \cite{MahalovTitiLeibovich90}. In that case, the increase of domain for the norms on the right-hand side is not needed.
\end{remark}

\begin{theo}\label{theo:uniquenesslocalenergyhelical}
If $u_0\in \mathring{E}_2$ is helical, the helical solution obtained in Theorem \ref{theo:globalweaklocalenergysolutions} is unique. Therefore, Theorem \ref{theo:globalexistenceintroduction} holds.
\end{theo}
\bpf 
We first claim that if $u$ is a helical local energy weak solution, then for any cylinder $D\subset \RRT$ $\pa_t u\in L^2((0,T),H^{-1}(D))$. Because $u$ is a weak solution, for any test function $\varphi\in C^\infty_c((0,T)\times D)$ we have 

$$
\int_0^T \int_D u \pa_t \varphi dxds = \int_0^T\int_D(u\otimes u+p-\na u)\na\varphi dx ds.$$
Therefore, we can bound as follows 
$$
\Big|\int_0^T \int_D u \pa_t \varphi dxds\Big| \lesssim \int_0^T (\norm{u}_{L^4_{uloc}}^2+\norm{\na u}_{L^2_{uloc}}+\norm{p}_{L^2_{uloc}})\norm \varphi_{H^1} ds.
$$
Using now Proposition \ref{prop:pressuredecomposition} and Lemma \ref{lem:helicalembedding}, we bound

$$
\Big|\int_0^T \int_D u \pa_t \varphi dxds\Big| \lesssim \int_0^T (\norm{u}_{L^2_{uloc}}\norm{\na u}_{L^2_{uloc}}+\norm{u}_{L^2_{uloc}}^2+\norm{\na u}_{L^2_{uloc}})\norm \varphi_{H^1} ds.
$$
We conclude that 
$$\Big|\int_0^T \int_D u \pa_t \varphi dxds\Big|\lesssim \norm{\varphi}_{L^2((0,T),H^1)}.$$
By duality, we prove the claim. 

We can now perform energy estimates on the difference $w=u-v$. The evolution equation is 
$$\pa_t w-\Delta w+\na\cdot (w\otimes u+v\otimes w)+\na(p-q)=0.$$
For any $k\in \RRT$, we define $\Xi_k=\chi_1(x-k)$. Multiply by $\Xi_k^2w$ and integrate in space. Integrating by parts and using that the solutions are divergence free, we obtain 
\begin{equation}\label{eq:uniquenessestimate}
\begin{aligned}
\frac{d}{dt}\int |\Xi_k w|^2+2\int |\Xi_k\na w|^2=&\int (\Delta \Xi_k^2)|w|^2-2\int \Xi_k^2 w \cdot (w\cdot\na)u+\int \na \Xi_k^2 \cdot v |w|^2 \\
&+ 2\int \na \Xi_k^2 \cdot w (p-q).
\end{aligned}
\end{equation}
Following Proposition \ref{prop:pressuredecomposition}, we get the bound 
$$\norm{p_{x_0}-q_{x_0}}_{L^q(C(x_0,1))}\lesssim \norm{w}_{L^{2q}_{\uloc}}(\norm{u}_{L^{2q}_{\uloc}}+\norm{v}_{L^{2q}_{\uloc}}).$$
We explain now how to deal with each term in the right hand side of \eqref{eq:uniquenessestimate}. The first one is directly bounded by $\norm{w}_{\LTu}$. In the second one, apply $L^4$ to each $\Xi_k w$, and $L^2_{\uloc}$ to $\na u$. Then use Lemma \ref{lem:helicalembedding}. Using the Young inequality for products, we have a term $\eps \norm{\Xi_k\na w}_{L^2}^2$ that can be absorbed on the LHS, and terms $(\norm{\Xi_k w}_{L}^2+\norm{w}_{\LTu}^2)\norm{\na u}_{\LTu}^2$. The third term is more troublesome. Applying Hölder and Lemma \ref{lem:helicalembedding} we do not get the weight $\Xi_k$ with $\na w$, all we get is $\norm{v}_{\LTu}\norm{w}_{\LTu}\norm{\na w}_{\LTu}$. For the last term, apply Hölder with exponents $4, 4/3$ and then use Proposition \ref{prop:pressuredecomposition} and Hölder with exponents $4,2$ to obtain the bound $\norm{w}_{L^4_{\uloc}}^2(\norm{u}_{\LTu}+\norm{v}_{\LTu})$. This allows us to get 
$$
\begin{aligned}
\frac{d}{dt}\int |\Xi_k w|^2+\int |\Xi_k\na w|^2\lesssim & 
\norm{\Xi_k w}_{L^2}^2\norm{\na u}_{\LTu}^2\\
&+\norm{w}_{\LTu}^2(\norm{\na u}_{\LTu}^2+1)+ \norm{w}_{\LTu}\norm{\na w}_{\LTu}(\norm{u}_{\LTu}+\norm{v}_{\LTu}).
\end{aligned}
$$
Note that the right-hand side is integrable by assumption. Denote $g(t)=\norm{\na u(t)}_{\LTu}^2$ and
$f(t)=\int_0^t\norm{w}_{\LTu}^2(\norm{\na u}_{\LTu}^2+1)+ \norm{w}_{\LTu}\norm{\na w}_{\LTu}(\norm{u}_{\LTu}+\norm{v}_{\LTu})ds$. Integrating the previous estimate, 

\begin{equation}\label{eq:uniquenessestimate}
\begin{aligned}
\norm{\Xi_k w(t)}_{L^2}^2 +\int_0^t \norm{\Xi_k\na w}^2_{L^2}ds\lesssim   & f(t)+\int_0^t \norm{\Xi_k w}_{L^2}^2g(s)ds.
\end{aligned}
\end{equation}
Denoting $y(t)=\int_0^t \norm{\Xi_k w}_{L^2}^2g(s)ds$, we observe that $y'(t)\lesssim g(t)f(t)+g(t)y(t)$. This ODE can be solved to obtain 
$$y(t)\lesssim e^{\int_0^tg(s)ds}f(t),$$
Then \eqref{eq:uniquenessestimate} and taking supremum in $k$ gives 
$$\norm{w}_{\LTu}^2(t) +\int_0^t \norm{\na w}^2_{\LTu}ds\lesssim e^{C\int_0^tg(s)ds}f(t).$$
Recall now that $f(t)=\int_0^t\norm{w}_{\LTu}^2(\norm{\na u}_{\LTu}^2+1)+ \norm{w}_{\LTu}\norm{\na w}_{\LTu}(\norm{u}_{\LTu}+\norm{v}_{\LTu})ds$. To treat the second term in $f$, apply the Young inequality for products with a small constant so that we can absorb the term with $\int_0^t\norm{\na w}_{\LTu}ds$ on the left-hand side. Then essentially we have 
$$\norm{w}_{\LTu}^2 +\int_0^t \norm{\na w}^2_{\LTu}\lesssim\int_0^th(s)\norm{w}_{\LTu}^2ds,$$
where $h$ is an integrable function. By Gronwall inequality, we conclude that $\norm{w(t)}_{\LTu}=0\quad\forall t\in [0,T].$
\epf

\newtheorem*{theo:globalwell-posednessintroduction}{Theorem \ref{theo:globalwell-posednessintroduction}}
\begin{theo:globalwell-posednessintroduction}[Global-in-time well-posedness for a helical vortex filament]
For any $\alpha\in\RR$ and any $\Gamma$ helical curve, there exists a unique global-in-time smooth helical solution to the Navier-Stokes equations with initial vorticity given by \eqref{vortexfilamentdata}.
\end{theo:globalwell-posednessintroduction}
\bpf
Due to Theorem \ref{theo:localexistence} and Corollary \ref{coro:uniquenesssmalltime}, we already know the local-in-time well-posedness. Because of Proposition \ref{prop:helicalsymmetrypreserved}, we know that our solution is helical. Note that the integrability of $u$ implies that $u(t_0)\in \mathring{E}_2$ for any $0<t_0$ small enough. Then, we can apply Theorem \ref{theo:globalexistenceintroduction} to uniquely extend the solution globally in time. 

For short times, our solution is smooth because it satisfies the Serrin condition for local regularity \cite[Chapter 13]{RobinsonRodrigoSadowski16}. Now, we do energy estimates to prove that the smoothness is preserved.

We keep the notation of the previous Theorem. We can find a time $t_1$ arbitrarily close to $t_0$ such that $\na u(t_1)\in \LTu(\RRT)$. We multiply the Navier-Stokes equations by $-\Xi_k^2\Delta u$ and integrate in space. We integrate by parts to bound $\frac{d}{dt}\norm{\Xi_k\na u}_{L^2}+\norm{\Xi_k\Delta u}_{L^2}$ by several terms. The general strategy to estimate terms is: Any $\Xi_k\Delta u$ is estimated in $L^2$, then apply the Young inequality for products to give $\norm{\Xi_k\Delta u}_{L^2}$ a small weight so that it can be absorbed on the left-hand side. For integrals with three terms, bound the lower-order ones in $L^4$ and apply Lemma \ref{lem:helicalembedding}. Due to the presence of the indicator function, other terms appear. To estimate the terms involving the pressure, pass the derivatives on the pressure to the other terms and use Proposition \ref{prop:pressuredecomposition}. To estimate the term with $\pa_t u$, substitute it using \eqref{NSvelocity} and proceed as before. The worst term remaining is

$$\frac{d}{dt}\norm{\Xi_k\na u}_{L^2}+\norm{\Xi_k\Delta u}_{L^2}\lesssim \norm{u}^2_{\LTu}\norm{\na u}^4_{\LTu},$$
which implies using Gronwall that $\na u(t)\in L^\infty([t_1,T];\LTu)$, $\Delta u(t)\in L^2([t_1,T];\LTu)$.

This is enough to apply the Serrin condition for local regularity to conclude that $u$ is smooth for all times.
\epf

\appendix
\section{Appendix}
\subsection{General analysis in $\RR^2\times\TT$}
\begin{lem}\label{lem:HolderBzLp}
	If $1\leq p_1,p_2,q\leq \infty$ and $\frac{1}{q}=\frac{1}{p_1}+\frac{1}{p_2}$, then 
	$$\norm{fg}_{B_{x_3}L^q_x}\lesssim \norm{f}_{B_zL^{p_1}_x}\norm{g}_{B_{x_3}L^{p_2}_x}.$$
\end{lem}
\bpf Write the expression of the left hand side. Apply Minkowsky's inequality in $x_3$ and Hölder inequality in $\widetilde{x}$.
\epf

\begin{lem}\label{lem:BSBzLp}
	$$\norm{\na\times(-\Delta)^{-1}f}_{B_{x_3}L^4_x}\lesssim \norm{f}_{B_{x_3}L^{4/3}_x}.$$
\end{lem}
\bpf Write the expression of the left hand side norm, then apply the change $\bar{x}\zeta=\widetilde{x}$. It reduces to proving standard estimates for $2D$ kernels.
\epf

\begin{lem}\label{lem:RieszBzLp}
	If $1< p< \infty$, for a constant depending on $p$
	$$\norm{(\Delta)^{-1}\na^2f}_{B_{x_3}L^p_x}\lesssim \norm{f}_{B_{x_3}L^p_x}.$$
\end{lem}
\bpf 
Fix a frequency $\zeta.$ For each $\zeta$ the inequality is true, so it also works for $B_{x_3}$.
\epf

\begin{lem}\label{lem:bztolp}
For $1<p\leq 2$, we have 
\begin{equation}\label{lembztolp}
\norm{f}_{B_{x_3}L^p_x}\lesssim \norm{f}_{L^p_{x}}^{1-\frac{1}{p}}\norm{\pa_{x_3}f}_{L^p_{x}}^{\frac{1}{p}} +\norm{f}_{L^p_{x}}
\end{equation}
\end{lem}
\bpf
We need the last term because $f$ might not have $0$ mean. Choose $N\in \NN$ and apply Hölder inequality in $\zeta$ to obtain
$$\sum_{|\zeta|=0}^N\norm{\hat{f}(\widetilde{x},\zeta)}_{L^p_x}\lesssim N^{\frac{1}{p}}\norm{\hat{f}(\widetilde{x},\zeta)}_{l^{p'}_\zeta L^p_x}\lesssim N^{\frac{1}{p}}\norm{\hat{f}(\widetilde{x},\zeta)}_{L^p_x l^{p'}_\zeta}\lesssim  N^{\frac{1}{p}}\norm{{f(x)}}_{L^p_{x}}.$$
Above, we can change the order of the norms because $p'\geq p$. We bound the high frequencies similarly by 
$$\sum_{|\zeta|>N}\norm{\hat{f}(\widetilde{x},\zeta)}_{L^p_x}\lesssim N^{\frac{1}{p}-1}\norm{\zeta\hat{f}(\widetilde{x},\zeta)}_{l^{p'}_\zeta L^p_x}\lesssim N^{\frac{1}{p}-1}\norm{{\pa_{x_3}f(x)}}_{L^p_{x}}.$$
Taking $N=\max\left\{1, \lfloor\frac{\norm{\pa_{x_3}f}_{L^p_{x}}}{\norm{f}_{L^p_{x,z}}}\rfloor \right\}$, we obtain the result.
\epf

\begin{lem}\label{lem:BiotSavartbounds}
Let $K(y)$ be the Biot-Savart kernel on $\RRT$. Then, for any multiindex $a\in\NN^3$,
$$
 \begin{aligned}
 |\pa_aK(y)|&\lesssim  |\widetilde{y}|^{-(1+a)}, \quad 1\leq |y| \\
 |\pa_aK(y)|&\lesssim  |y|^{-(2+a)}, \quad |y|\leq 1.
 \end{aligned}
 $$
\end{lem}
\bpf
A different proof can be seen in \cite[Lemma 3.3]{bronzinussenzveig2015}. We show that the kernel in $\RRT$ is obtained summing in the third component the 3D kernel, i.e.,
\begin{equation}\label{eq:biotsavartsum}
K(\widetilde{y},y_3)=\frac{1}{4\pi}\sum_{k\in\ZZ}\frac{(y_1,y_2,y_3+2k\pi)}{(|\widetilde{y}|^2+|y_3+2\pi k|^2)^{3/2}}.
\end{equation}
This sum is $y_3$-periodic. Passing the vector norm inside the sum, we bound $$|K(\widetilde{y},y_3)|\lesssim\sum_{k\in\ZZ}\frac{1}{|\widetilde{y}|^2+|y_3+2\pi k|^2}.$$
Since $y_3\in [-\pi,\pi)$, we have that $|y_3+2\pi k|$ is always equal or greater than $|y_3|$. Then we can bound
$$\sum_{k\in\ZZ}\frac{1}{|\widetilde{y}|^2+|y_3+2\pi k|^2}\lesssim \sum_{k=0}^\infty\frac{1}{|\widetilde{y}|^2+|y_3+2\pi k|^2}\lesssim \frac{1}{|\widetilde{y}|^2+|y_3|^2}+\int_1^\infty \frac{dy_3}{|\widetilde{y}|^2+|y_3|^2}\lesssim \frac{1}{|\widetilde{y}|+1}+\frac{1}{|\widetilde{y}|^2+|y_3|^2}.$$
This shows that the sum \eqref{eq:biotsavartsum} is absolutely convergent away from $0$, so the definition makes sense. We need to show that such $K$ solves $\curl K=\delta_{0}(1,1,1)$. We know that this is true for the term $k=0$ of the sum, while for the other terms the result is identically $0$. Due to absolute convergence, we can switch sum and derivatives. The bound for derivatives of $K$ using this sum is analogue.
\epf

\begin{lem}\label{lem:Heatproductspace}
Let $H_D$ be the heat kernel on the space $D$, $\widetilde{y}\in\RR^2$, $y_3\in\TT$. Then 
$$\Heat_\RRT(t,\widetilde{y},y_3)=\sum_{k\in\ZZ}\Heat_{\RR^3}(t,\widetilde{y},y_3+2\pi k)=\Heat_{\RR^2}(t,\widetilde{y})\sum_{k\in\ZZ}\Heat_\RR(t,y_3+2\pi k)=\Heat_{\RR^2}(t,\widetilde{y})\Heat_\TT(t,y_3).$$
\end{lem}
\bpf To prove the first equality we can proceed as in the previous Lemma. We can separate the variables due to properties of the exponential function. 
\epf

\subsection{Analysis close to the curve $\Gam$}

\begin{prop}\label{prop:estBSeta}
    Let $\eta$ be supported in $\{|\widetilde{x}|\leq 16r\}$. If $0<r\ll 1$ is small enough, then the following holds:
\begin{equation}\label{estBSeta}
\norm{\chi_{8r}P_\Phi^{-1}(-\Delta)^{-1}\na\times Q_\Phi\eta}_{B_{x_3}L^4_x}+\norm{\na(\chi_{r}P_\Phi^{-1}(-\Delta)^{-1}\na\times Q_\Phi\eta)}_{B_{x_3}L^{4/3}_x}\lesssim \norm{\eta}_{B_{x_3}L^{4/3}_x}
\end{equation}
\begin{equation}\label{estBSawayeta}
\norm{(1-\chi_r)\chi_{8r}P_\Phi^{-1}(-\Delta)^{-1}\na\times Q_\Phi(\chi_{\frac{r}{4}}\eta)}_{B_{x_3}L^4_x}\lesssim r^{-\frac{1}{2}}\norm{\eta}_{B_{x_3}L^1_x}
\end{equation}
\end{prop}
\begin{prop}\label{prop:estdiffBS}
Let $m>1$, $0<r\ll 1$ be small enough, $\eta$ supported in $\{|\widetilde{x}|\leq 16r\}$. Assuming $0<\sqrt{t}\lesssim r\ll 1$ and 
$$v=\chi_rP_\Phi^{-1}\na\times(-\Delta)^{-1}Q_\Phi\eta,$$
we have 
\begin{equation}\label{estdiffBS}
    t^{\frac{1}{4}}\norm{v-\na\times(-\Delta)^{-1}\eta}_{B_{x_3}L^4_x}+
    t^{\frac{1}{4}}\norm{\na(v-\na\times(-\Delta)^{-1}\eta)}_{B_{x_3}L^{4/3}_x}\lesssim (t^\frac{1}{4}r^{-\frac{1}{2}}+r\ln{r})\norm{\eta}_{N^0}
\end{equation}
\end{prop}

\begin{lem}\label{lem:boundfourierK}
Let $K$ be the Biot-Savart kernel in $\RR^2\times\TT.$ Then the corresponding kernel $\widetilde{K}$ of the operator $\FF Q_\Phi^{-1}\chitil_{16r}KQ_\Phi\chi_{16r}\FF$ satisfies 
$$|\widetilde{K}(\widetilde{x},\zeta,\widetilde{x}',\zeta')|\lesssim |\widetilde{x}-\widetilde{x}'|^{-1}\< \zeta-\zeta'\>^{-2}.$$
\end{lem}

\bpf
The three previous results are proven in \cite[Appendix A]{bedrossiangermainharropgriffiths23} for $x\in \RR^2\times\TT$, $y\in\RR^3$. We have instead $y\in \RR^2\times\TT$, so Biot-Savart law and the inverse of the Laplacian change in this space. We can still use the same proof because we can restrict $|y|<K$ due to the presence of $\chi_{cr},\chitil_{cr}$ through the computations, so essentially the kernels behave identically. 
\epf

\begin{lem}\label{lem:boundduhamelheat} It holds
\begin{enumerate} 
    \item For $f_1$ with support in $\{|\widetilde{x}|\leq 16r\}$, 
$$
\begin{aligned}
&\norm{Q_\Phi^{-1}\left(\chitil_{8r}\int_0^te^{(t-s)\Delta}(Q_\Phi \divergence f_1(s)+\divergence f_2(s)\right)ds}_{N^0} \\
&\hspace{3cm} \lesssim \sup_{0<s\leq t} \left( s^\frac{3}{4}\norm{\< s^{-\frac{1}{2}}\widetilde{x}\>^m f_1(s)}_{B_{x_3}L_x^{4/3}}+s^\frac{3}{4}\norm{\< s^{-\frac{1}{2}}\bd \>^m f_2(s)}_{L_y^{2}} \right).
\end{aligned}
$$    
\item For $f$ with support in $\{|x|\leq 4r\}$ and any multi-index $\gamma \in \NN^3$ with $|\gamma | \leq 2$,

$$
\norm{Q_\Phi^{-1}\left(\chitil_{8r}\int_0^te^{(t-s)\Delta}Q_\Phi \na_{x,z}^\gamma f(s)ds\right)}_{N^0} \lesssim t^{1-\frac{|\gamma|}{2}}\norm{f}_{N^\beta}.
$$ 
\end{enumerate}
\end{lem}
\bpf
See \cite[Lemma 7.12]{bedrossiangermainharropgriffiths23}. In this case $f_2$ can be nonzero far from $\Gam$, but the result is still valid since only $L^p$ norms of the heat kernel are used, and in both $\RR^3$ and $\RRT$ we have $\norm{\na^k e^{t\Delta}}_{L^p}\lesssim t^{-\frac{k}{2}-\frac{3(p-1)}{2p}}$.
\epf

\subsection{Morrey spaces}
Here we proof some results for the spaces defined in Definition \ref{def:morreyR2T}.

\begin{lem}\label{lem:MorreyHolder}
If $\frac{1}{p}=\frac{1}{p_1}+\frac{1}{p_2}$ and $\frac{1}{q}=\frac{1}{q_1}+\frac{1}{q_2}$
$$ 
\norm{fg}_{M^p_q}\leq\norm{f}_{M^{p_1}_{q_1}}\norm{g}_{M^{p_2}_{q_2}}.
$$
\end{lem}
\bpf
Use the usual Hölder inequality.
\epf

\begin{lem}\label{lem:Morreyembeddings}
$$\begin{aligned}   
L^p&\subset M^p_q \quad \forall \, 1\leq q \leq p,\\
B_{x_3}L^q_x&\subset M_q^{\frac{3q}{2}}.
\end{aligned}$$
\end{lem}
\bpf
The first embedding is a consequence of the Hölder inequality. For the second one, apply Hölder in the $x_3$ variable and use $B_{x_3}\subset L^\infty$.
\epf

\subsection*{Acknowledgments.} 		
F. Gancedo and A. Hidalgo-Torné were partially supported by the ERC through the Starting Grant H2020-EU.1.1.-639227, by the Fundacion de Investigación de la Universidad de Sevilla through the grant FIUS23/0207, by the MICINN (Spain) through the grants EUR2020-112271 and PID2020-114703GB-I00 and by the Junta de Andalucía through the grant P20-00566. F. Gancedo was partially supported by MINECO grant RED2018-102650-T (Spain). F. Gancedo acknowledges support
from IMAG, funded by MICINN through the Maria de Maeztu
Excellence Grant CEX2020-001105-M/AEI/10.13039/501100011033. A. Hidalgo-Torné was partially supported by the grant FPU19/02748, funded by the Spanish Ministry of Universities. 

\bibliographystyle{plain}
\bibliography{biblioteca}

\end{document}